\documentclass[12pt]{article}
\UseRawInputEncoding
\usepackage{mathrsfs}
\usepackage{amsmath}
\usepackage{amssymb}
\usepackage{amsthm}
\usepackage{amsfonts}
\usepackage{amscd}
\usepackage[mathscr]{eucal}
\newcommand{\Z} {{\mathbb  Z}}
\newcommand{\Q}{{\mathbb  Q}}

\begin{document}
\parindent  25pt
\baselineskip  10mm
\textwidth  15cm    \textheight  23cm \evensidemargin -0.06cm
\oddsidemargin -0.01cm

\title{ {On algebraic congruence varieties over semirings }}
\author{\mbox{}
{ Derong Qiu }
\thanks{ \quad E-mail:
derong@mail.cnu.edu.cn, \ derongqiu@gmail.com } \\
(School of Mathematical Sciences,
 Capital Normal University, \\
Beijing 100048, P.R.China )  }

\date{}
\maketitle
\parindent  24pt
\baselineskip  10mm
\parskip  0pt

\par   \vskip 0.4cm

{\bf Abstract} \quad In this paper, we develop some foundations
for a theory of algebraic varieties of congruences on commutative semirings.
By studying the structure of congruences, firstly, we show that
the spectrum $ \text{Spec}^{c}(A) $ consisting of prime congruences
on a semirings $ A $ has a Zariski topological structure; Then, for
two semirings $ A \subset B, $ we consider the polynomial semiring
$ S = A[x_{1}, \cdots , x_{n}] $ and the affine $ n-$space $ B^{n}. $
For any congruence $ \sigma $ on $ S $ and congruence $ \rho $ on $ B, $
we introduce the $ \rho-$algebraic varieties
$ Z_{\rho }(\sigma )(B)$ in $ B^{n}, $ which are the set of zeros
in $ B^{n} $ of the system of polynomial $ \rho-$congruence equations
given by $ \sigma . $ When $ \rho $ is a prime congruence, we find these
varieties satisfying the axiom of closed sets, and forming a (Zariski) topology
on $ B^{n}. $ Some results about their structures
including a version of Nullstellensatz of congruences are obtained.
\par  \vskip  0.4 cm

\par  \vskip  0.4 cm


\par     \vskip  0.6 cm

\hspace{-0.8cm}{\bf 1. Introduction}

\par \vskip 0.2 cm

As known, ideals are very important in studying ring structure and the
polynomial equations over rings, and, in a commutative ring, the ideal
and the congruence corresponds to each other equivalently. Yet, in a
commutative semiring, this changes, and congruences seem to behave more
natural than ideals. \\
Our motivation here is to consider the problem of solving polynomial
equations over semirings, like the same important problem in
algebraic geometry on rings. We focus on congruences, and try to develop
some foundations for a theory of algebraic varieties of congruences on
commutative semirings. Firstly, we study prime congruence,
maximal congruence, congruence generated by a subset and a few other types
congruence, and obtain some results about their structures
(see Prop.2.6, 2.10, 2.12, 2.13, 2.16, 2.19, 2.21, 2.23).
In particular, we show that the spectrum $ \text{Spec}^{c}(A) $
consisting of prime congruences on a semirings $ A $ has a Zariski
topological structure (see Theorem 2.25). Then, for
two semirings $ A \subset B, $ we consider the polynomial semiring
$ S = A[x_{1}, \cdots , x_{n}] $ and the affine $ n-$space $ B^{n}. $
For any congruence $ \sigma $ on $ S $ and congruence $ \rho $ on $ B, $
we introduce the $ \rho-$algebraic varieties
$ Z_{\rho }(\sigma )(B)$ in $ B^{n}, $ which are the set of zeros
in $ B^{n} $ of the system of polynomial $ \rho-$congruence equations
given by $ \sigma $ (see Def.3.1). When $ \rho $ is a prime congruence,
we find these varieties satisfying the axiom of closed sets,
and forming a (Zariski) topology on $ B^{n} $ (see Theorem 3.4).
Some results about their structures including a version of
Nullstellensatz of congruences are obtained (see Prop.3.9,
Theorems 3.7, 3.10, 3.11, 3.12). \\
These facts show one possible way among others to consider
the problem of solving polynomial (congruence) equations in semirings,
they seem to be well corresponded to the ones of classical algebraic geometry
of ideals on commutative rings, yet in a different background.
Also, the results of (Zariski) topology for congruences and varieties
established here, might be useful in further study on some more deep questions
about such as Sheaf theory, scheme on semirings. \\
Another possible application of this work, is about some questions in tropical
geometry. One of the key connection might be the Nullstellensatz question.
Recently, in the tropical setting, an analogue of Hilbert's Nullstellensatz
has been extensively studied and many interesting results about it
have been obtained (see, e.g., [BE], [JM1], [Iz], [IzS], [GG], [Gr], [MR], etc.).
In this paper, we provide a version of Nullstellensatz of $ \rho-$algebraic
varieties on semirings. For an possible application in tropical geometry,
it seems that the radical of congruence used in this paper need to be improved,
and we wish to consider it in a future work. \\
There are more general version of algebraic equations than our congruence
ones considered here. When we do some work on arithmetic problems in number
theory, we usually only need to consider the solution of positive integers
even prime numbers. So I ask myself a question that how to work on algebraic
equations if there is no subtraction? and I come to realize that the
usual polynomial equations can be viewed as a special class of the general
algebraic relations, which led to my consideration of the present work.
In fact, I learn that this is related to the context of universal algebraic
geometry and logical geometry (see, e.g., [P1-2],[PP],[DMR1-2]). \\
Lastly, we need to give an explanation about the use of prime congruences
here. Unlike working with prime ideals on rings, when we study congruences
on semirings, we confront with several versions of prime congruences
in the literature, and they are usually used for different purposes.
The prime congruences of Def.2.1(1) below we use in this work based on
our trying to find a suitable topology to do some geometry about
the set of (congruence) zeros, i.e., the varieties, as shown in
Theorem 3.4 below. It seems to us that not all the other version of
prime congruences can satisfy our demand here. \\
{\bf Notation.} \ The sign $ \subset $ is used for inclusion of a subset,
including the possibility of equality. When we say that $ A $ is a proper
subset of $ B, A \subsetneqq B $ is intended.

\par     \vskip  0.4 cm

\hspace{-0.8cm}{\bf 2. The structure of congruences }

\par \vskip 0.2 cm

Let $ (A, +, \cdot ) $ be a semiring, that is, $ A $ is a non-empty set
on which we have defined operations of addition and multiplication satisfying
the following four conditions (see [G, p.1]): \
(1) \ $ (A, +) $ is a commutative monoid with identity element $ 0; $ \
(2) \ $ (A, \cdot ) $ is a monoid with identity element $ 1 \neq 0; $ \
(3) \ $ a(b + c) = ab + ac $ and $ (a + b)c = ac + bc \
(\forall a, b ,c \in A); $ \
(4) \ $ 0a = 0 = a0 \ (\forall a \in A). $ \\
The semiring $ A $ is commutative if the monoid $ (A, \cdot ) $ is
commutative. In the following, unless otherwise stated, the semiring
always means a commutative semiring. \\
Now for a commutative semiring $ A, $ recall that a ( binary) relation
on $ A $ is a subset of $ A \times A $ (see [Ho, p.14]), and a relation
$ \rho $ of $ A $ is a congruence if it is an equivalence relation
satisfying the following condition (see [G, p.4]):  \\
$ (a, b), (c, d) \in \rho \Longrightarrow
(a + c, b + d) \in \rho  $ and $ (ac, bd) \in \rho \ (a, b, c, d \in A). $ \\
So an equivalence relation $ \rho $ is a congruence on $ A $ if and only if
$ \rho $ is a subsemiring of $ A \times A. $ The identity congruence on $ A $
is denoted by $ \text{id}_{A} = \{(a, a) : a \in A \}. $ Let $ \rho $ be a
congruence on $ A, $ if $ \rho \neq A \times A $ then it is called a proper
congruence. Obviously, $ \rho $ is not proper if and only if $ (1, 0) \in \rho. $
The quotient semiring of $ A $ by $ \rho $ is denoted by
$ \overline{A} = A / \rho . $ For $ c \in A, $ we denote $ \overline{c} =
c \ \text{mod} \rho \in \overline{A}. $ Moreover, for any
positive integer $ n $ and subset $ C $ of $ A^{n}, $ where $ A^{n} =
\{(a_{1}, \cdots , a_{n}) : \ a_{1}, \cdots , a_{n} \in A \}, $ we denote
$ C / \rho = \{(\overline{c_{1}}, \cdots , \overline{c_{n}}) : \
(c_{1}, \cdots , c_{n})  \in C \} \subset (A / \rho )^{n}. $
For $ 0 \neq a \in A, $ if there exists an element $ 0 \neq b \in A $
such that $ a \cdot b = 0, $ then $ a $ is called a zero divisor of $ A. $
If $ A $ contains no zero divisors, then $ A $ is called a semidomain.
If $ A $ is a semidomain, and $ (A \setminus \{0 \}, \cdot ) $ is a
multiplicative group, then $ A $ is called a semifield. For example,
$ (\Z_{\geq 0}, +, \cdot ) $ is a semidomain,
and $ (\Q_{\geq 0}, +, \cdot ) $ is a semifield, where $ \Z_{\geq 0} $
(resp. $ \Q_{\geq 0}$) is the set of non-negative integers
(resp. rational numbers).
\par  \vskip 0.2cm

{\bf Definition 2.1.} \ Let $ A $ be a commutative semiring, and $ \rho $
be a congruence of $ A $ such that $ \rho \neq A \times A. $ \\
(1) \ If $ \rho $ satisfies the following condition: \\
for $ (a, b), (c, d) \in A \times A, \ (ac + bd, ad + bc) \in \rho
\Rightarrow (a, b) \in \rho $ or $ (c, d) \in \rho. $ \\
Then $ \rho $ is called a prime congruence. \\
(2) \ If $ \rho $ satisfies the following condition: \\
$ (ab, 0) \in \rho  \Rightarrow (a, 0) \in \rho $ or $ (b, 0) \in \rho
\quad (a, b \in A). $ \\
Then $ \rho $ is called a semi-prime congruence. \\
(3) \ If $ \rho $ satisfies the following condition: \\
for a congruence $ \tau: \ \rho \subset \tau \subset A \times A \Rightarrow
\tau = \rho $ or $ \tau = A \times A. $ \\
Then $ \rho $ is called a maximal congruence. \\
(4) \ If $ \rho $ satisfies the following condition: \\
for $ a \in A, \ (a, 0) \notin \rho  \Rightarrow (ab, 1) \in \rho $ for
some $ b \in A. $ \\
Then $ \rho $ is called a semi-maximal congruence. \\
Note that there are several other definitions for prime congruences
(yet their meaning are not always the same) in the literature on semirings.
The above Def.2.1(1) of prime congruences is the same as in [JM1-2].
We denote
\begin{align*} &\text{Spec}^{c}(A) = \{\rho : \ \rho \ \text{is a prime
congruence of} \ A \}, \\
&\text{Spec}^{c^{\prime }}(A) =
\{\rho : \ \rho \ \text{is a semi-prime congruence of} \ A \}, \\
&\text{Max}^{c}(A) = \{\rho : \ \rho \ \text{is a maximal congruence
of} \ A \}, \\
&\text{Max}^{c^{\prime }}(A) = \{\rho : \ \rho \ \text{is a semi-maximal
congruence of} \ A \}.
\end{align*}
And call $ \text{Spec}^{c}(A) $ (resp. $ \text{Spec}^{c^{\prime }}(A),
\text{Max}^{c}(A) $ or $ \text{Max}^{c^{\prime }}(A)$) the spectrum of
prime (resp. semi-prime, maximal or semi-maximal) congruences of $ A. $
Obviously, $ \text{Spec}^{c}(A) \subset \text{Spec}^{c^{\prime }}(A). $
\par  \vskip 0.2cm

Recall that an ideal $ I $ of a commutative semiring $ A $
is a non-empty subset $ I $ of $ A $ which is closed under addition and
satisfies the condition that if $ a \in A $ and $ b \in I $ then $ ab \in I. $
A prime ideal $ P $ of $ A $ is an ideal $ P $ which satisfies the condition
that if $ a \cdot b \in P $ then $ a \in P $ or $ b \in P \ (a, b \in A). $
We denote $ \mathcal{I}(A) = \{\text{all ideals of} \ A \}, \
\mathcal{E}(A) =  \{\text{all equivalences of} \ A \}$ and
$ \mathcal{C}(A) = \{\text{all congruences of} \ A \}. $
For any non-empty subsets $ S $  and $ T $ of $ A, $ we denote
$ S + T = \{s + t : \ s \in S, \ t \in T \}, \ S \cdot T =
\{s t : \ s \in S, \ t \in T \}, $ and for $ a \in A,
a + S = \{a + s : \ s \in S \}, \  a \cdot S = \{a s : \ s \in S \}. $
For $ J \in \mathcal{I}(A) $ and $ \sigma \in \mathcal{C}(A), $ we denote
$ \rho _{J} = \{(a, b) \in A \times A : \ a + J = b + J \} $ and
$ I_{\sigma } = \{a \in A : \ (a, 0) \in \sigma \}. $
\par  \vskip 0.2cm

{\bf Lemma 2.2.} \ Let $ A $ be a commutative semiring. We have \\
(1) \ $ \rho _{J} \in \mathcal{C}(A) \ (\forall J \in \mathcal{I}(A)); $
\quad (2) \ $ I_{\sigma } \in \mathcal{I}(A) \ (\forall \sigma \in
\mathcal{C}(A)); $ \\
(3) \ $ I_{\rho _{J}} \subset J \ (\forall J \in \mathcal{I}(A)); $
\quad (4) \ $ \rho _{I_{\sigma }} \subset \sigma \ (\forall \sigma \in
\mathcal{C}(A)); $ \\
(5) \ $ J \subset J^{\prime } \Rightarrow \rho _{J} \subset
\rho _{J^{\prime }} \ (J, J^{\prime } \in \mathcal{I}(A)); $
\quad (6) \ $ \sigma \subset \tau \Rightarrow I_{\sigma }
\subset I_{\tau } \ (\sigma , \tau \in \mathcal{C}(A)). $
\par  \vskip 0.1cm
{\bf Proof.} \ Follow directly from the definitions. \quad $ \Box $ \\
Note that if $ A $ is a commutative ring (obviously it is also a
commutative semiring), then it is easy to see that
$ I_{\rho _{J}} = J $ and $ \rho _{I_{\sigma }} = \sigma $ for all
$ J \in \mathcal{I}(A) $ and $ \sigma \in \mathcal{C}(A). $ But this
property changes for a commutative semiring, i.e., the inclusions in
the above Lemma 2.2(3-4) may be strict in general for a commutative
semiring $ A. $ For example, let $ A = (\Z_{\geq 0}, +, \cdot) $ be
the semidomain as above. Take the ideal
$ J = 2 \cdot \Z_{\geq 0} = \{2 n : \ n \in \Z_{\geq 0} \}. $ Then
it is easy to see that $ \rho _{J} = \text{id}_{ A}. $
Hence $ I_{\rho _{J}} = \{0 \} \subsetneqq J, $ i.e.,
$ I_{\rho _{J}} \subsetneqq J. $
Now take the congruence $ \sigma = \{(a, b) \in A \times A : \ a \equiv
b \ (\text{mod} \ 2) \} \neq \text{id}_{ A}. $ Then by definition,
$ I_{\sigma } = \{a \in A : \ a \equiv 0 \ (\text{mod} \ 2) \} =
J. $ So $ \rho _{I_{\sigma }} = \rho _{J} = \text{id}_{ A} \subsetneqq \sigma , $
i.e., $ \rho _{I_{\sigma }} \subsetneqq \sigma . $ \\
From now on, we mainly study the congruences on commutative semirings. \\
For a congruence $ \rho \neq A \times A $ on a commutative semiring
$ A, $ it is easy to see that, \\
$ \rho $ is a semi-prime congruence $ \Leftrightarrow A / \rho $ is a semidomain; \\
$\rho $ is a semi-maximal congruence $ \Leftrightarrow A / \rho $ is a semifield; \\
$ \rho $ is contained in a maximal congruence of $ A. $ In particular,
$ A $ contains a maximal congruence. Obviously, a semi-maximal congruence is
also a semi-prime congruence. Moreover, it is easy to see that the intersection
of a family of congruences of $ A $ is also a congruence. So, for any
relation $ R $ on $ A, $ there exists a unique smallest congruence of $ A $
containing $ R, $ which is the intersection of all congruences of $ A $
containing $ R. $ We denote it by $ R^{c}, $ and call it the congruence
generated by $ R. $
\par  \vskip 0.2cm

{\bf Definition 2.3.} \ Let $ A $ be a commutative semiring. We
use the operation $ \ast $ to denote the twisted product on
$ A \times A $ as follows (see e.g.,[BE], [JM1-2]):
$$ (a, b) \ast (c, d) = (ac + bd, ad + bc) \quad (\forall (a, b), (c, d)
\in A \times A). $$
Then for every positive integer $ n, $ we define $ (a, b)^{\ast n} $
inductively as follows: \\
$ (a, b)^{\ast 1} = (a, b), \ (a, b)^{\ast 2} = (a, b) \ast (a, b), \
\cdots , (a, b)^{\ast n} = (a, b)^{\ast (n - 1)} \ast (a, b). $ \\
We also define $ (a, b)^{\ast 0} = (1, 0). $ \\
Using this operation, the condition for a congruence $ \rho $ of $ A $
to be a prime congruence is as follows: \
$ (a, b) \ast (c, d) \in \rho \Rightarrow (a, b) \in \rho $ or $ (c, d)
\in \rho \ (a, b, c, d \in A). $ \\
In the following, for two non-empty subsets $ S $ and $ T $ of
$ A \times A, $ we shall write $ S \ast T = \{\alpha \in A \times A : \
\alpha = (a, b) \ast (c, d) \ \text{with} \ (a, b) \in S,
(c, d) \in T \}. $
\par  \vskip 0.2cm

{\bf Lemma 2.4.} \ Let $ A $ be a commutative semiring. \\
(1) \ $ \ast $ is associative: \ $ ((a, b) \ast (c, d)) \ast (e, f)
= (a, b) \ast ((c, d) \ast (e, f)); $  \\
(2) \ $ \ast $ is commutative: \ $ (a, b) \ast (c, d) =
(c, d) \ast (a, b); $ \\
(3) \ $ \ast $ is distributive: \ $ ((a, b) + (c, d)) \ast (e, f) =
(a, b) \ast (e, f) + (c, d) \ast (e, f); $ \\
(4) \ $ (a, b) \ast (c, c) \in \text{id}_{A}; $ \\
(5) \ $ (a, b) \ast (1, 0) = (a, b), \ (a, b) \ast (0, 1) = (b, a), \
(a, b) \ast (0, 0) = (0, 0). $ \\
(6) \ for positive integers $ n, $ we have \\
$ (b, a)^{\ast n} =
\left \{\begin{array}{l} (a, b)^{\ast n} \quad \quad \text{if} \
2 \mid n  \\
(a, b)^{\ast n} \ast (0, 1) \quad \text{if} \
2 \nmid n.
\end{array} \right. $ In particular,
$ (0, 1)^{\ast n} =
\left \{\begin{array}{l} (1, 0) \quad \text{if} \
2 \mid n  \\
(0, 1) \quad \text{if} \ 2 \nmid n.
\end{array} \right. $ \\
$ (a, b, c, d, e, f \in A). $
\par  \vskip 0.1cm
{\bf Proof.} \ (1)$\sim$ (5) follow directly from the definitions, and
(6) is proved by induction. \quad $ \Box $
\par  \vskip 0.2cm

{\bf Remark 2.5.} \ Let $ A $ be a commutative semiring. Then $ (A \times A, +, \ast ) $
is a commutative semiring with additive zero
$ (0, 0) $ and multiplicative unity $ (1, 0), $ and the
congruences on $ A $ may connect with the ideals of $ (A \times A, +, \ast ). $
\par  \vskip 0.2cm

{\bf Proposition 2.6.} \ Let $ A $ be a commutative semiring, and
$ n $ be a positive integer. Then for $ (a, b), (c, d) \in A \times A, $
we have
\begin{align*} &(a, b)^{\ast n} = (\sum _{2 \mid i, \ i = 0}^{n}
\begin{pmatrix}  n \\
i \end{pmatrix} a^{n - i}
 b^{i}, \ \sum _{2 \nmid i, \ i = 1}^{n} \begin{pmatrix}  n \\
i \end{pmatrix} a^{n - i} b^{i}), \\
&((a, b) + (c, d))^{\ast n} = (a + c, b + d)^{\ast n} =
\sum _{i = 0}^{n} \begin{pmatrix}  n \\
i \end{pmatrix} (a, b)^{\ast i} \ast (c, d)^{\ast (n - i)}.
\end{align*}
{\bf Proof.} \ We use induction on $ n. $ The case for $ n = 1 $ is obvious.
For the first equality, assume it holds for $ n, $ then
\begin{align*} &(a, b)^{\ast (n + 1)} = (a, b) \ast (a, b)^{\ast n} =
(a, b) \ast (\sum _{2 \mid i, \ i = 0}^{n}
\begin{pmatrix}  n \\
i \end{pmatrix} a^{n - i} b^{i}, \ \sum _{2 \nmid i, \ i = 1}^{n}
\begin{pmatrix}  n \\
i \end{pmatrix} a^{n - i} b^{i}) \\
&=(\sum _{2 \mid i, \ i = 0}^{n} \begin{pmatrix}  n \\
i \end{pmatrix} a^{n + 1 - i} b^{i} + \sum _{2 \nmid i, \ i = 1}^{n}
\begin{pmatrix}  n \\
i \end{pmatrix} a^{n - i} b^{i + 1}, \ \sum _{2 \mid i, \ i = 0}^{n}
\begin{pmatrix}  n \\
i \end{pmatrix} a^{n - i} b^{i + 1} \\
& \quad \quad \quad \quad + \sum _{2 \nmid i, \ i = 1}^{n} \begin{pmatrix}  n \\
i \end{pmatrix} a^{n + 1 - i} b^{i}) \\
&= (\sum _{2 \mid i, \ i = 0}^{n} \begin{pmatrix}  n + 1 \\
i \end{pmatrix} a^{n + 1 - i} b^{i} + \frac{1 + (-1)^{n + 1}}{2}b^{n + 1},
\ \sum _{2 \nmid i, \ i = 1}^{n} \begin{pmatrix}  n + 1 \\
i \end{pmatrix} a^{n + 1 - i} b^{i} \\
& \quad \quad \quad \quad + \frac{1 - (-1)^{n + 1}}{2}b^{n + 1}) \\
&= (\sum _{2 \mid i, \ i = 0}^{n + 1} \begin{pmatrix}  n + 1 \\
i \end{pmatrix} a^{n + 1 - i} b^{i}, \ \sum _{2 \nmid i, \ i = 1}^{n + 1}
\begin{pmatrix}  n + 1 \\
i \end{pmatrix} a^{n + 1 - i} b^{i}).
\end{align*}
So via the induction, the first equality holds for all positive
integers $ n. $ \\
For the second equality, assume it holds for $ n, $ then
\begin{align*} &((a, b) + (c, d))^{\ast (n + 1)}= ((a, b) + (c, d)) \ast
\sum _{i = 0}^{n} \begin{pmatrix} n \\
i \end{pmatrix} (a, b)^{\ast i} \ast (c, d)^{\ast (n - i)} \\
&= \sum _{i = 0}^{n} \begin{pmatrix} n \\
i \end{pmatrix} (a, b)^{\ast (i + 1)} \ast (c, d)^{\ast (n - i)} +
\sum _{i = 0}^{n} \begin{pmatrix} n \\
i \end{pmatrix} (a, b)^{\ast i} \ast (c, d)^{\ast (n + 1 - i)}  \\
&= \sum _{i = 1}^{n + 1} \begin{pmatrix} n \\
i - 1 \end{pmatrix} (a, b)^{\ast i} \ast (c, d)^{\ast (n + 1 - i)} +
\sum _{i = 0}^{n} \begin{pmatrix} n \\
i \end{pmatrix} (a, b)^{\ast i} \ast (c, d)^{\ast (n + 1 - i)}  \\
&= \sum _{i = 1}^{n} \begin{pmatrix} n + 1 \\
i \end{pmatrix} (a, b)^{\ast i} \ast (c, d)^{\ast (n + 1 - i)} +
\begin{pmatrix} n + 1 \\
n + 1 \end{pmatrix} (a, b)^{\ast (n + 1)} \ast (c, d)^{\ast 0} \\
&+ \begin{pmatrix} n + 1 \\
0 \end{pmatrix} (a, b)^{\ast 0} \ast (c, d)^{\ast (n + 1)} \\
&= \sum _{i = 0}^{n + 1} \begin{pmatrix}  n + 1 \\
i \end{pmatrix} (a, b)^{\ast i} \ast (c, d)^{\ast (n + 1 - i)}.
\end{align*}
So via the induction, the second equality holds for all positive
integers $ n. $ \quad $ \Box $
\par  \vskip 0.2cm

{\bf Definition 2.7.} \ Let $ A $ be a commutative semiring and $ R $ be
a relation of $ A. $ We define \ $ R_{+} = \{(a, b) \in A \times A : \
(a + c, b + c) \in R \ \text{for some} \ c \in A\}. $
\par  \vskip 0.2cm

{\bf Lemma 2.8.} \ Let $ \rho $ be a congruence on a commutative semiring
$ A. $ Then $ \rho _{+} $ is also a congruence of $ A, $
and $ \rho _{+} \supset \rho. $ Moreover, $ (\rho _{+})_{+} = \rho _{+}. $
\par  \vskip 0.1cm
{\bf Proof.} \ Follows easily from the definition.  \quad $ \Box $
\par  \vskip 0.2cm

{\bf Definition 2.9.} \ Let $ \rho $ be a congruence on a commutative semiring
$ A. $ We define a radical of $ \rho $ by
$$ \sqrt{\rho } = \{(a, b) \in A \times A : \ (a + c, b + c)^{\ast n}
\in \rho \ \text{for some} \ c \in A \
\text{and some positive integer} \ n \}. $$
This can also be defined as follows, as pointed out
by an anonymous expert:
$$ \sqrt{\rho } = \{(a, b) \in A \times A : \ (a, b)^{\ast n} + (c, c)
\in \rho \ \text{for some} \ c \in A \
\text{and some positive integer} \ n \}, $$
which follows from the fact that
$ (a + c + 1, b + c +  1)^{\ast n} = (a, b)^{\ast n} + (c, c) + (d, d) $
for some $ d \in A $ by the binomial expansion in Prop.2.6 above
(note that $ (c + 1, c +  1)^{\ast n} = (c, c) +
((2^{n-1} n - 1)c, (2^{n-1} n - 1)c ) +
(e, e) $ for some $ e \in A $). Obviously, $ \rho \subset \sqrt{\rho }. $
\par  \vskip 0.2cm

{\bf Proposition 2.10.} \ Let $ \rho $ be a congruence on a commutative
semiring $ A. $ Then $ \sqrt{\rho } $ is also a congruence on $ A. $
Moreover, $ (\sqrt{\rho })_{+} = \sqrt{\rho }. $
\par  \vskip 0.1cm
{\bf Proof.} \ Firstly, $ (a, b) \in \sqrt{\rho } \Leftrightarrow
(b, a) \in \sqrt{\rho }. $ In fact, if $ (a, b) \in \sqrt{\rho }, $
then $ (a + c, b + c)^{\ast n} \in \rho $ for some $ c \in A $ and some
positive integer $ n. $ So $ (a + c, b + c)^{\ast n} \ast (0, 1) \in \rho . $
Thus by Lemma 2.4 above, $ (b + c, a + c)^{\ast n} \in \rho , $ so
$ (b, a) \in \sqrt{\rho }. $ and vice versa. Next, if $ (a, b), (b, c) \in
\sqrt{\rho }, $ then $ (a + e, b + e)^{\ast m} \in \rho $ and
$ (b + f, c + f)^{\ast n} \in \rho $ for some $ e, f \in A $ and some positive
integers $ m, n. $ Then by Lemma 2.4 and Prop.2.6 above,
\begin{align*} &(a + e + f, b + e + f)^{\ast m} =
((a + e, b + e) + (f, f))^{\ast m} \\
&= \sum _{i = 0}^{m} \begin{pmatrix} m \\
i \end{pmatrix} (a + e, b + e)^{\ast i} \ast (f, f)^{\ast (m - i)} \\
&= (a + e, b + e)^{\ast m} + \sum _{i = 0}^{m - 1} \begin{pmatrix} m \\
i \end{pmatrix} (a + e, b + e)^{\ast i} \ast (f, f)^{\ast (m - i)} \in \rho .
\end{align*}
Similarly, $ (b + e + f, c + e + f)^{\ast n} \in \rho . $ Write
$ t = e + f, $ then
$ (a + t, b + t)^{\ast m}, \ (b + t, c + t)^{\ast n} \in \rho . $
Write $ s = b + 2t \in A $ and $ k = m + n. $ For $ i = 0, 1, \cdots , k, $
if $ i \geq m, $ then $ (a + t, b + t)^{\ast i}
= (a + t, b + t)^{\ast m} \ast (a + t, b + t)^{\ast (i - m)} \in \rho . $
If $ i < m, $ then $ k - i > k - m = n, $ so
$ (b + t, c + t)^{\ast (k - i)} = (b + t, c + t)^{\ast n} \ast
(b + t, c + t)^{\ast (k - i - n)} \in \rho . $ Then
\begin{align*} &(a + s, c + s)^{\ast k} =
((a + t, b + t) + (b + t, c + t))^{\ast k} \\
&= \sum _{i = 0}^{k} \begin{pmatrix} k \\
i \end{pmatrix} (a + t, b + t)^{\ast i} \ast (b + t, c + t)^{\ast (k - i)}
\in \rho .
\end{align*}
So $ (a, c) \in \sqrt{\rho }. $ Also $ \text{id}_{A} \subset \rho \subset
\sqrt{\rho }. $ Therefore, $ \sqrt{\rho } $ is an equivalence relation of $ A. $ \\
Next, let $ (a, b), (c, d) \in \sqrt{\rho }. $ Then by definition,
$ (a + e, b + e)^{\ast m} \in \rho $ and $ (c + f, d + f)^{\ast n}
\in \rho $ for some $ e, f \in A $ and some positive integers $ m, n. $ Write
$ k = m + n. $ By Prop.2.6 above, $ ((a + e, b + e) + (c + f, d + f))^{\ast k} =
\sum _{i = 0}^{k} \begin{pmatrix}  k \\
i \end{pmatrix} (a + e, b + e)^{\ast i} \ast (c + f, d + f)^{\ast (k - i)}. $
For $ i = 0, 1, \cdots, k, $ if $ i \geq m, $ then $ (a + e, b + e)^{\ast i}
= (a + e, b + e)^{\ast m} \ast (a + e, b + e)^{\ast (i - m)} \in \rho ; $
if $ i < m, $ then $ k - i > k - m = n, $ so $ (c + f, d + f)^{\ast (k - i)}
= (c + f, d + f)^{\ast n} \ast (c + f, d + f)^{\ast (k - i - n)} \in \rho , $
which implies that $ ((a + e, b + e) + (c + f, d + f))^{\ast k} \in \rho , $
i.e., $ ((a + c + e + f, b + d + e + f))^{\ast k} \in \rho . $ So
$ (a + c, b + d) \in \sqrt{\rho }, $ i.e., $ (a, b) + (c, d) \in \sqrt{\rho }. $
The remainder is to show that $ (a, b) \in \sqrt{\rho } \Rightarrow (ac, bc)
\in \sqrt{\rho } \ (\forall c \in A). $ Indeed, $ (a, b) \in \sqrt{\rho }
\Rightarrow (a + e, b + e)^{\ast n} \in \rho $ for some $ e \in A $ and some
positive integer $ n. $ Then by Prop.2.6 above,
\begin{align*} &((a + e)c, (b + e)c)^{\ast n} \\
&= (\sum _{2 \mid i, \ i = 0}^{n}
\begin{pmatrix} n \\
i \end{pmatrix} ((a + e)c)^{n - i} \cdot ((b + e)c)^{i}, \
\sum _{2 \nmid i, \ i = 1}^{n} \begin{pmatrix} n \\
i \end{pmatrix} ((a + e)c)^{n - i} \cdot ((b + e)c)^{i}) \\
&= c^{n} \cdot (\sum _{2 \mid i, \ i = 0}^{n}
\begin{pmatrix} n \\
i \end{pmatrix} (a + e)^{n - i} \cdot (b + e)^{i}, \
\sum _{2 \nmid i, \ i = 1}^{n} \begin{pmatrix} n \\
i \end{pmatrix} (a + e)^{n - i} \cdot (b + e)^{i}) \\
&= c^{n} \cdot (a + e, b + e)^{\ast n} \in \rho  \
(\text{here we write} \ x \cdot (y, z) = (xy, xz) \ (\forall x, y, z \in A)),
\end{align*}
so $ (ac+ ec, bc + ec)^{\ast n} \in \rho , $ hence $ (ac, bc) \in \sqrt{\rho }. $
Therefore, $ \sqrt{\rho } $ is a congruence. To show that $ (\sqrt{\rho })_{+} =
\sqrt{\rho }, $ firstly, by Lemma 2.8 above, $ \sqrt{\rho } \subset
(\sqrt{\rho })_{+}. $ Conversely, if $ (a, b) \in (\sqrt{\rho })_{+}, $ then
$ (a + c, b + c) \in \sqrt{\rho } $ for some $ c \in A. $ So
$ (a + c + d, b + c + d)^{\ast n} \in \rho $ for some $ d \in A $ and some positive
integer $ n, $ i.e., $ (a + e, b + e)^{\ast n} \in \rho $ with $ e = c + d \in A, $
so $ (a, b) \in \sqrt{\rho }, $ and so $ (\sqrt{\rho })_{+} \subset \sqrt{\rho }. $
Therefore $ (\sqrt{\rho })_{+} = \sqrt{\rho }, $ and the proof is completed.
\quad $ \Box $
\par  \vskip 0.2cm

{\bf Definition 2.11.} \ Let $ A $ and $ \rho $ be as in Def.2.9 above. \\
(1) \ If $ \sqrt{\rho } = \rho , $ then $ \rho $ is called
to be a radical congruence on $ A. $ \\
(2) \ If $ \sqrt{\rho } = \rho _{+}, $ then $ \rho $ is called
to be a quasi-radical congruence on $ A. $ \\
(3) \ Denote $ R_{\text{nil}}(A) = \{(a, b) \in A \times A : \ (a, b)^{\ast n}
\in \text{id}_{A} \ \text{for some positive integer} \ n \} $ and
$ \rho _{\text{nil}}(A) = R_{\text{nil}}(A)_{+}, $ then
$ \rho _{\text{nil}}(A) = \sqrt{\text{id}_{A}}, $ and $ \rho _{\text{nil}}(A) $
is called to be the nilpotent congruence of $ A. $
\par  \vskip 0.2cm

{\bf Proposition 2.12.} \ Let $ A $ be a commutative semiring. Let $ \rho ,
\rho _{1} $ and $ \rho _{2} $ be congruences on $ A. $ \\
(1) \ $ \rho _{1} \subset \rho _{2} \Rightarrow \sqrt{\rho _{1}} \subset
\sqrt{\rho _{2}} $ and $ (\rho _{1})_{+} \subset (\rho _{2})_{+}. $ \quad
(2) \ $ \rho _{+} \subset \sqrt{\rho }. $ \\
(3) \ $ \sqrt{\sqrt{\rho }} = \sqrt{\rho } = \sqrt{\rho _{+}}. $
Particularly, $ \sqrt{\rho } $ is a radical congruence on $ A. $  \\
(4) \ $ (a, b), \ (b, c) \in R_{\text{nil}}(A) \ \Rightarrow \ (a, c)^{\ast k} +
(e, e) \in \text{id}_{A} $ for some $ e \in A $ and some positive
integer $ k. $ \\
(5) \ If $ \rho $ is a prime congruence, then
$ \sqrt{\rho } = \rho _{+}, $ i.e., $ \rho $ is a quasi-radical congruence.
If the prime congruence $ \rho $ satisfies $ \rho = \rho _{+}, $ then
$ \rho $ is a radical congruence.
\par  \vskip 0.1cm
{\bf Proof.} \ (1) and (2) follow directly from the definitions. \\
(3) \ Since $ \rho \subset \sqrt{\rho }, $ we only need to show that
$ \sqrt{\sqrt{\rho }} \subset \sqrt{\rho }. $ Let
$ (a, b) \in \sqrt{\sqrt{\rho }}, $ then $ (a, b)^{\ast n} + (c, c) \in
\sqrt{\rho } $ for some $ c \in A $ and some $ n \in \Z_{> 0}. $ So
$ ((a, b)^{\ast n} + (c, c))^{\ast m} + (d, d) \in \rho $ for some
$ d \in A $ and some $ m \in \Z_{> 0}. $ By the binomial expansion in
Prop.2.6 above, we get $ ((a, b)^{\ast n} + (c, c))^{\ast m} =
(a, b)^{\ast (mn)} + (e, e) $ with $ e \in A, $ which implies
$ (a, b) \in \sqrt{\rho }, $ so $ \sqrt{\sqrt{\rho }} = \sqrt{\rho }. $
Also $ \rho \subset \rho _{+} \subset \sqrt{\rho }, $
so $ \sqrt{\rho } \subset \sqrt{\rho _{+}} \subset \sqrt{\sqrt{\rho }}, $
and so $\sqrt{\rho _{+}} =  \sqrt{\rho }. $ \\
(4) \ Follows from the above Prop.2.6 and Lemma 2.4. \\
(5) \ Follows from the definition. The proof is completed. \quad $ \Box $ \\
Recall that the semiring $ A $ satisfies the additive annihilation law if \
$ a + c = b + c \Rightarrow a = b \ (a, b, c \in A). $ \
For example, $ (\Z_{\geq 0}, +, \cdot ) $ satisfies the additive annihilation
law.
\par  \vskip 0.2cm

{\bf Proposition 2.13.} \ Let $ A $ be a commutative semiring satisfying the
additive annihilation law. Then $ \rho _{\text{nil}}(A) = R_{\text{nil}}(A), $
and $ \rho _{\text{nil}} (A / \rho _{\text{nil}}(A)) = \text{id}. $
\par  \vskip 0.1cm
{\bf Proof.} \ Follows from the above Prop.2.6 and Prop.2.10. \quad $ \Box $ \\
For a commutative semiring $ A, $ if $ R_{\text{nil}}(A) = \text{id}_{A}, $
then $ A $ is called to be reduced. If $ \rho _{\text{nil}}(A) = \text{id}_{A}, $
then $ A $ is called to be strongly reduced. obviously, strongly reduced
$ \Rightarrow $ reduced. We denote $ \mathfrak{N}_{c}(A) = R_{\text{nil}}(A)^{c}, $
which is the congruence of $ A $ generated by the relation $ R_{\text{nil}}(A), $
and call $ \mathfrak{N}_{c}(A) $ the quasi-nilpotent congruence of $ A. $
Obviously, $  \rho _{\text{nil}}(A) \supset  \mathfrak{N}_{c}(A), $ also,
$ A $ is reduced if and only if $ \mathfrak{N}_{c}(A) = \text{id}_{A}. $ \\
Now for two relations $ R $ and $ R^{\prime } $ on a non-empty set $ S, $
recall that their product $$ R \circ R^{\prime } = \{(a, c) \in S \times S : \
\exists \ b \in S \ \text{such that} \ (a, b) \in
R \ \text{and} \ (b, c) \in R^{\prime } \}, $$ and the inverse $ R^{-1} =
\{(a, b) \in S \times S : \ (b, a) \in R \} $ (see [Ho, pp.14, 15]).
Obviously, for relations $ R, R _{1}, \cdots, R_{n}, R^{\prime },
R^{\prime \prime } $ on $ S, $ one has $ (R^{-1})^{-1} = R; \
(R_{1} \circ \cdots \circ R_{n})^{-1} = R_{n}^{-1} \circ \cdots \circ R_{1}^{-1}; \
R \subset R^{\prime } \Rightarrow R^{-1} \subset R^{\prime ^{-1}}, \
R \circ R^{\prime \prime } \subset R^{\prime } \circ R^{\prime \prime } $
and $ R^{\prime \prime } \circ R \subset R^{\prime \prime } \circ R^{\prime }. $ \\
As in [Ho], for any relation $ R $ on $ S, $ we denote its transitive
closure by $ R^{\infty }, $ which is defined by $ R^{\infty } =
\cup _{n = 1}^{\infty } R^{n}, $ where $ R^{n} =
\underbrace{R \circ \cdots \circ R}_{n}. $ It is well known that
$  R^{\infty } $ is the smallest transitive relation on $ S $ containing
$ R $ (see [Ho, p.20]). Also, we denote by $ R^{e} $ the equivalence
on $ S $ generated by $ R, $ i.e., the smallest equivalence on $ S $
containing $ R. $ Then $ R^{e} = (R \cup R^{-1} \cup \text{id}_{S})^{\infty } $
(see [Ho, p.20]). Moreover, $ (x, y) \in R^{e} \Leftrightarrow x = y $ or
for some positive integer $ n $ there exists a sequence $ x = z_{1}
\rightarrow z_{2} \rightarrow \cdots \rightarrow z_{n} = y $ in which for
each $ i \in \{1, 2, \cdots , n - 1\} $ either $ (z_{i}, z_{i + 1})
\in R $ or $ (z_{i + 1}, z_{i}) \in R $ (see [Ho, p.21]).
\par  \vskip 0.2cm

{\bf Definition 2.14.} \ Let $ A $ be a commutative semiring, $ R $ be
a relation on $ A, $ and $ E $ be an equivalence relation on $ A. $
We define the relations
\begin{align*} & R^{L} =
\{(ax + y, bx + y) : \ (a, b) \in R, x, y \in A \}, \ \text{and} \\
&E^{\flat} = \{(a, b) \in A \times A : \ (ax + y, bx + y) \in E \
\text{for all} \ x, y \in A \}.
\end{align*}
Since $ 0, 1 \in A, $ we have $ R \subset R^{L}. $ Also, $ \rho ^{L}
= \rho $ for any congruence $ \rho $ of $ A. $
\par  \vskip 0.2cm

{\bf Lemma 2.15.} \ Let $ A $ be a commutative semiring, $ R_{1} $
and $ R_{2} $ be two relations on $ A, $ and $ E $ be an equivalence
relation on $ A. $ We have \\
(1) \ $ R_{1} \subset R_{1}^{L}; $ \quad (2) \ $ (R_{1}^{L})^{-1} =
 (R_{1}^{-1})^{L}; $ \quad (3) \ $ R_{1} \subset R_{2} \Rightarrow
R_{1}^{L} \subset R_{2}^{L}; $ \\
(4) \ $ (R_{1}^{L})^{L} = R_{1}^{L}; $ \quad (5) \ $ (R_{1} \cup R_{2})^{L}
= R_{1}^{L} \cup R_{2}^{L};  $ \\
(6) \ $ R_{1} = R_{1}^{L} $ if and only if $ R_{1} $ satisfies the
following condition: \\
$ (a, b) \in R_{1} \Rightarrow (ax, bx), (a + y, b + y) \in R_{1} \
(\forall x, y \in A). $ \\
(7) \ $ R_{1} = R_{1}^{L} \Rightarrow R_{1}^{n} = (R_{1}^{n})^{L} $
for all positive integer $ n. $ \\
(8) \ $ E^{\flat} $ is the largest congruence on $ A $ contained in $ E. $
\par  \vskip 0.1cm
{\bf Proof.} \ Follow easily from the definitions. \quad $ \Box $
\par  \vskip 0.2cm

{\bf Proposition 2.16.} \ Let $ R $ be a relation on a commutative semiring
$ A. $ Then $ R^{c} = (R^{L})^{e}, $ where $ R^{c} $ is the congruence of
$ A $ generated by $ R $ as before.
\par  \vskip 0.1cm
{\bf Proof.} \ By definition, $ (R^{L})^{e} $ is the equivalence on $ A $
generated by $ R^{L}. $ Denote $ S = R^{L} \cup (R^{L})^{-1} \cup \text{id}_{A}, $
then $ (R^{L})^{e} = S^{\infty } $ (see [Ho, p.20]). By Lemma 2.15 above,
$ S = (R \cup R^{-1} \cup \text{id}_{A})^{L}, $ so $ S^{L} = S, $ and so
$ S^{n} = (S^{n})^{L} $ for all positive integers $ n. $ Now for any
$ (a, b) \in (R^{L})^{e}, $ by the above discussion, $ (a, b) \in S^{n} $
for some positive integer $ n. $ Since $ S^{n} = (S^{n})^{L}, $ by
Lemma 2.15 above, $ (ax, bx) \in S^{n} $ and $ (a + y, b + y) \in S^{n}, $
hence $ (ax, bx), (a + y, b + y) \in (R^{L})^{e} \ (\forall x, y \in A). $
In particular, if $ (a, b), (c, d) \in (R^{L})^{e}, $ then $ (a + c, b + c)
\in (R^{L})^{e} $ and $ (b + c, b + d) \in (R^{L})^{e}, $ so
$ (a + c, b + d) \in (R^{L})^{e}. $ Hence $ (R^{L})^{e} $ is a congruence
of $ A $ containing $ R. $ Now let $ \rho $ be a congruence of $ A $
containing $ R. $ Then by Lemma 2.15 above, $ R^{L} \subset \rho ^{L} = \rho . $
So $ (R^{L})^{e} \subset \rho ^{e} = \rho , $ which shows that
$ R^{c} = (R^{L})^{e}, $ and the proof is completed. \quad $ \Box $
\par  \vskip 0.2cm

{\bf Corollary 2.17.} \ Let $ R $ be a relation on a commutative semiring
$ A. $ Then for $ a, b \in A, (a, b) \in R^{c} $ if and only if $ a = b $ or
for some positive integer $ n $ there exists a sequence $ a = z_{1}
\rightarrow z_{2} \rightarrow \cdots \rightarrow z_{n} = b $ in which for
each $ i \in \{1, 2, \cdots , n - 1\} $ either $ (z_{i}, z_{i + 1})
\in R^{L} $ or $ (z_{i + 1}, z_{i}) \in R^{L}. $
\par  \vskip 0.1cm
{\bf Proof.} \ Follows from the above Prop.2.16 (see [Ho, p.21]).
\quad $ \Box $ \\
Let $ R $ and $ \rho $ be a relation and a congruence on a commutative
semiring $ A, $ respectively. Recall that \
$ R/ \rho = \{(a\rho , b\rho) : \ (a, b) \in R \} = \{(\overline{a},
\overline{b}) \in A/ \rho \times A/ \rho : \ (a, b) \in R \} $ \ is a
relation on the quotient semiring $ \overline{A} = A/ \rho . $ If
$ \text{id}_{A} \subset R \subset \rho , $ then
$ R/ \rho = \text{id}_{\overline{A}}. $
\par  \vskip 0.2cm

{\bf Corollary 2.18.} \ Let $ R $ and $ \rho $ be a relation and a
congruence on a commutative semiring $ A, $ respectively. Then
$ R^{c} / \rho = (R/ \rho )^{c} $ on $ \overline{A} = A/ \rho . $
\par  \vskip 0.1cm
{\bf Proof.} \ Follows from the above Cor.2.17.  \quad $ \Box $  \\
For a commutative semiring $ A, $ recall that $ \mathcal{E}(A) $ and
$ \mathcal{C}(A) $ denote the set of all equivalences and the set of
all congruences on $ A, $ respectively. It is known that both
$ (\mathcal{E}(A), \subset, \wedge, \vee) $ and
$ (\mathcal{C}(A), \subset, \wedge, \vee) $ are complete lattices
(for lattice, see e.g.,[Ho, p.12]): if $ \rho , \sigma
\in \mathcal{E}(A), $ then $ \rho \wedge \sigma = \rho \cap \sigma $
is their intersection, while $ \rho \vee \sigma = (\rho \cup \sigma )^{e}. $
Similarly, for $ \rho , \sigma \in  \mathcal{C}(A), \rho \wedge \sigma =
\rho \cap \sigma $ and $ \rho \vee \sigma = (\rho \cup \sigma )^{c} =
((\rho \cup \sigma )^{L})^{e} = (\rho \cup \sigma )^{e} $ follows by
Lemma 2.15 and Prop.2.16 above. So the join
of $ \rho $ and $ \sigma $ in $ \mathcal{C}(A) $ coincides
with their join in $ \mathcal{E}(A). $ Both $ \mathcal{E}(A) $ and
$ \mathcal{C}(A) $ have the maximum element $ A \times A $ and minimum
element $ \text{id}_{A}. $
\par  \vskip 0.2cm

{\bf Proposition 2.19.} \ Let $ \rho $ and $ \sigma $ be two congruences on
a commutative semiring $ A. $ Then $ \rho \vee \sigma =
(\rho \circ \sigma )^{\infty }. $ In other words, for $ a, b \in A, \ (a, b) \in
\rho \vee \sigma $ if and only if for some positive integer $ n $ there
exist elements $ x_{1}, x_{2}, \cdots , x_{2n - 1} \in A $ such that
$ (a, x_{1}) \in \rho , (x_{1}, x_{2}) \in \sigma , (x_{2}, x_{3}) \in \rho ,
\cdots , (x_{2n - 2}, x_{2n -1}) \in \rho, (x_{2n - 1}, b) \in \sigma . $
\par  \vskip 0.1cm
{\bf Proof.} \ Follows from Prop.2.16 above.  \quad $ \Box $
\par  \vskip 0.2cm

{\bf Corollary 2.20.} \ If $ \rho $ and $ \sigma $ are two congruences on
a commutative semiring $ A $ satisfying $ \rho \circ \sigma =
\sigma \circ \rho . $ Then $ \rho \vee \sigma = \rho \circ \sigma . $
\par  \vskip 0.1cm
{\bf Proof.} \ Follows from Prop.2.19 above.  \quad $ \Box $ \\
For example, if $ A $ is a semifield and $ \rho, \sigma \in \mathcal{C}(A), $
then $ \rho \circ \sigma = \sigma \circ \rho, $ so
$ \rho \vee \sigma = \rho \circ \sigma. $
\par  \vskip 0.2cm

{\bf Proposition 2.21.} \ Let $ A $ be a commutative semiring. For $ a, b \in A, $
denote \\
$ R(a, b) = \{(ax + by + z, bx + ay + z) : \ x, y, z \in A \} =
(a, b) \ast (A \times A) + \text{id}_{A}. $ \\
Then $ R(a, b)_{+} $ (as defined in Def.2.7 above) is a congruence on $ A, $
and $ R(a, b) \subset \rho \subset R(a, b)_{+} \subset \rho _{+} $ with
$ \rho = \{(a, b) \}^{c}. $
\par  \vskip 0.1cm
{\bf Proof.} \ Obviously, $ \text{id}_{A} \subset R(a, b)_{+}. $ Also,
if $ (u, v) \in R(a, b)_{+}, $ then $ (u + s, v + s) \in R(a, b) $
for some $ s \in A, $ so $ (u + s, v + s) = (ax + by + z, bx + ay + z) $
for some $ x, y, z \in A, $ so $ (v + s, u + s) = (ay + bx + z, by + ax + z)
\in R(a, b), $ so $ (v, u) \in R(a, b)_{+}. $ Now let $ (u, v), (v, w)
\in R(a, b)_{+}. $ Then $ (u + s, v + s), (v + t, w + t) \in R(a, b) $
for some $ s, t \in A, $ so $ (u + s, v + s) = (ax + by + z, bx + ay + z) $
and $ (v + t, w + t) = (ax^{\prime } + by^{\prime } + z^{\prime },
bx^{\prime } + ay^{\prime } + z^{\prime }) $ for some $ x, y, z, x^{\prime },
y^{\prime }, z^{\prime } \in A. $ Then $ (u + s + t, v + s + t) =
(ax + by + z + t, bx + ay + z + t), $ and $ (v + s + t, w + s + t) =
(ax^{\prime } + by^{\prime } + z^{\prime } + s,
bx^{\prime } + ay^{\prime } + z^{\prime } + s). $ Write $ d = v + s + t =
bx + ay + z + t = ax^{\prime } + by^{\prime } + z^{\prime } + s, $ then
$ u + s + t + d = a (x + x^{\prime }) + b (y + y^{\prime }) + z + z^{\prime }
+ s + t, $ and $ w + s + t + d = b (x + x^{\prime }) + a (y + y^{\prime }) +
z + z^{\prime } + s + t. $ So $ (u + s + t + d, w + s + t + d) \in R(a, b), $
and so $ (u, w) \in R(a, b)_{+}. $ Therefore, $ R(a, b)_{+} $ is an equivalence
on $ A. $ Next, let $ (u, v), (u^{\prime }, v^{\prime }) \in R(a, b)_{+}, $
then $ (u + s, v + s), (u^{\prime } + s^{\prime }, v^{\prime } + s^{\prime })
\in R(a, b) $ for some $ s, s^{\prime } \in A. $ So $ (u + s, v + s) =
(ax + by + z, bx + ay + z) $ and $ (u^{\prime } + s^{\prime },
v^{\prime } + s^{\prime }) = (ax^{\prime } + by^{\prime } + z^{\prime },
bx^{\prime } + ay^{\prime } + z^{\prime }) $ for some $ x, y, z, x^{\prime },
y^{\prime }, z^{\prime } \in A. $ Hence $ (u + u^{\prime } + s + s^{\prime },
v + v^{\prime } + s + s^{\prime }) = (a(x + x^{\prime }) + b(y + y^{\prime })
+ z + z^{\prime }, b(x + x^{\prime }) + a(y + y^{\prime }) + z + z^{\prime })
\in R(a, b), $ so $ (u, v) + (u^{\prime }, v^{\prime }) =
(u + u^{\prime }, v + v^{\prime }) \in R(a, b)_{+}. $ Moreover, for
$ (u, v), s, x, y, z $ above, let $ t \in A, $ then $ (ut + st, vt + st) =
(a(tx) + b(ty) + tz, b(tx) + a(ty) + tz) \in R(a, b), $ so $ (ut, vt)
\in R(a, b)_{+}. $ Therefore, $ R(a, b)_{+} $ is a congruence on $ A. $
The inclusion $ R(a, b) \subset \rho $ follows from Cor.2.5 above. The
remainder can be verified directly, and the proof is completed.
\quad $ \Box $
\par  \vskip 0.2cm

{\bf Lemma 2.22.} \ Let $ A $ be a commutative semiring and $ \rho $ be a
congruence on $ A. $ let $ (a, b), (c, d) \in A \times A $ and $ c = z_{1}
\rightarrow z_{2} \rightarrow \cdots \rightarrow z_{n} = d \ (n \geq 2) $
be a sequence in $ A. $ If $ (a, b) \ast (z_{i}, z_{i + 1}) \in \rho _{+} $
for all $ i \in \{1, 2, \cdots , n - 1 \}, $ then
$ (a, b) \ast (c, d) \in \rho _{+}. $
\par  \vskip 0.1cm
{\bf Proof.} \ By Lemma 2.8 above, $ \rho _{+} $ is a congruence on $ A, $
and $ (\rho _{+})_{+} = \rho _{+}. $ Write $ t = a z_{2} + b z_{2}, $
then from $ (a, b) \ast (z_{1}, z_{2}) \in \rho _{+} $ and
$ (a, b) \ast (z_{2}, z_{3}) \in \rho _{+}, $ we have
$ (az_{1} + bz_{3} + t, az_{3} + bz_{1} + t)
= (a, b) \ast (z_{1}, z_{2}) + (a, b) \ast (z_{2}, z_{3}) \in \rho _{+}, $
so $ (az_{1} + bz_{3}, az_{3} + bz_{1}) \in (\rho _{+})_{+} = \rho _{+}, $
i.e., $ (a, b) \ast (z_{1}, z_{3}) \in \rho _{+}. $ Then, from
$ (a, b) \ast (z_{1}, z_{3}) \in \rho _{+} $ and $ (a, b) \ast (z_{3}, z_{4})
\in \rho _{+}, $ we also have
$ (a, b) \ast (z_{1}, z_{4}) \in \rho _{+}, $ and so on. After finite
steps, we get $ (a, b) \ast (z_{1}, z_{n}) \in \rho _{+}, $ i.e.,
$ (a, b) \ast (c, d) \in \rho _{+}. $ The proof is completed.
\quad $ \Box $
\par  \vskip 0.2cm

{\bf Proposition 2.23.} \ Let $ A $ be a commutative semiring and $ \rho $ be a
maximal congruence on $ A. $ If $ \rho = \rho _{+}, $ then $ \rho $ is
also a prime congruence.
\par  \vskip 0.1cm
{\bf Proof.} \ Let $ (a, b), (c, d) \in A \times A $ be such that
$ (a, b) \ast (c, d) \in \rho . $ If $ (c, d) \notin \rho , $ then we
need to show that $ (a, b) \in \rho . $ To see this, let $ R =
\{(c, d) \} \cup \rho , $ then $ R $ is a relation on $ A $ and
$ R \supsetneqq \rho . $ So $ R^{c} \supsetneqq \rho , $ and so
$ R^{c} = A \times A $ because $ \rho $ is maximal. In particular,
$ (1, 0) \in R^{c}. $ By Cor.2.17 above, there exists a sequence $ 1 = z_{1}
\rightarrow z_{2} \rightarrow \cdots \rightarrow z_{n} = 0 $ in which for
each $ i \in \{1, 2, \cdots , n - 1 \} $ either $ (z_{i}, z_{i + 1})
\in R^{L} $ or $ (z_{i + 1}, z_{i}) \in R^{L}. $ By Lemma 2.15 above,
$ R^{L} = (\{(c, d) \} \cup \rho )^{L} = \{(c, d) \}^{L} \cup \rho ^{L} =
\{(c, d) \}^{L} \cup \rho $ because $ \rho ^{L} = \rho . $ By definition,
$ \{(c, d) \}^{L} = \{(cx + y, dx + y) : \ x, y \in A \}. $ So
$ R^{L} = \{(cx + y, dx + y) : \ x, y \in A \} \cup \rho . $ As above, for
each $ i \in \{1, 2, \cdots , n - 1 \} $ either $ (z_{i}, z_{i + 1})
\in R^{L} $ or $ (z_{i + 1}, z_{i}) \in R^{L}. $ If $ (z_{i}, z_{i + 1})
\in R^{L}, $ then $ (z_{i}, z_{i + 1}) \in \rho $ or $ (z_{i}, z_{i + 1})
= (cx_{i} + y_{i}, dx_{i} + y_{i}) $ for some $ x_{i}, y_{i} \in A. $
If $ (z_{i}, z_{i + 1}) \in \rho , $ then $ (a, b) \ast (z_{i}, z_{i + 1})
\in \rho . $ If $ (z_{i}, z_{i + 1}) = (cx_{i} + y_{i}, dx_{i} + y_{i}), $
write $ t = ay_{i} + by_{i} \in A, $ then \
$ (a, b) \ast (z_{i}, z_{i + 1}) = (a, b) \ast (cx_{i} + y_{i}, dx_{i} + y_{i})
= (acx_{i} + bdx_{i} + t, adx_{i} + bcx_{i} + t)
= ((a, b) \ast (c, d)) \cdot (x_{i}, x_{i}) + (t, t) \in \rho. $ \
So $ (z_{i}, z_{i + 1}) \in R^{L} \Rightarrow (a, b) \ast (z_{i}, z_{i + 1})
\in \rho . $ \\
Similarly, $ (z_{i + 1}, z_{i}) \in R^{L} \Rightarrow (a, b) \ast (z_{i + 1}, z_{i})
\in \rho . $ \\
Note that $ (a, b) \ast (z_{i + 1}, z_{i}) = ((a, b) \ast (z_{i}, z_{i + 1}))
\ast (0, 1). $ So the above discussion shows that
$ (a, b) \ast (z_{i}, z_{i + 1}) \in \rho $ for all
$ i \in \{1, 2, \cdots , n - 1 \}. $ By our assumption, $ \rho = \rho _{+}. $
Hence by Lemma 2.22 above, we get
$ (a, b) \ast (1, 0) = (a, b) \ast (z_{1}, z_{n}) \in \rho , $ i.e.,
$ (a, b) \in \rho . $ Therefore, $ \rho $ is a prime congruence on $ A, $
and the proof is completed. \quad $ \Box $
\par  \vskip 0.2cm

{\bf A question.} Is a maximal congruence in a commutative semiring always
prime?
\par  \vskip 0.2cm

As known, in [JM1] it is shown that maximal
congruences are always prime in an additively idempotent semiring.
\par  \vskip 0.2cm

{\bf Definition 2.24.} \ Let $ A $ be a commutative semiring and $ \sigma $
be a congruence on $ A. $ We define $ V^{\text{co}}(\sigma ) =
\{\rho \in \text{Spec}^{c}(A) : \ \rho \supset \sigma \}. $
\par  \vskip 0.2cm

{\bf Theorem 2.25.} \ Let $ A $ be a commutative semiring. \\
(1) \ If $ \sigma _{1} $ and $ \sigma _{2} $ are two congruences on $ A, $
then $ V^{\text{co}}(\sigma _{1}) \cup  V^{\text{co}}(\sigma _{2}) =
 V^{\text{co}}(\sigma _{1} \cap \sigma _{2}). $ \\
(2) \ If $ \{\sigma _{\alpha } \}_{\alpha \in \Lambda } $ is a family of
congruences on $ A, $ then
$ \cap _{\alpha \in \Lambda } V^{\text{co}}(\sigma _{\alpha }) =
V^{\text{co}}(\sigma ), $ where $ \sigma =
(\cup _{\alpha \in \Lambda } \sigma _{\alpha })^{c} $
is the congruence on $ A $ generated by the set
$ \cup _{\alpha \in \Lambda } \sigma _{\alpha }. $ \\
(3) \ $ V^{\text{co}}(\text{id}_{A}) = \text{Spec}^{c}(A) $ and
$ V^{\text{co}}(A \times A) = \emptyset . $
\par  \vskip 0.1cm
{\bf Proof.} \ (1) \ Since $ \sigma _{1} \cap \sigma _{2} \subset
\sigma _{i} \ (i = 1, 2), \ V^{\text{co}}(\sigma _{1} \cap \sigma _{2})
\supset V^{\text{co}}(\sigma _{i}) \ (i = 1, 2), $ so
$ V^{\text{co}}(\sigma _{1} \cap \sigma _{2}) \supset
V^{\text{co}}(\sigma _{1}) \cup  V^{\text{co}}(\sigma _{2}). $ Conversely,
if $ \rho \in V^{\text{co}}(\sigma _{1} \cap \sigma _{2}), $ then
$ \rho \supset \sigma _{1} \cap \sigma _{2}. $ If $ \rho \nsupseteq
\sigma _{1}, $ equivalently, $ \rho \notin V^{\text{co}}(\sigma _{1}), $
then there exists an element $ (a, b) \in \sigma _{1}, $ but
$ (a, b) \notin \rho . $ Now for any $ (c, d) \in \sigma _{2}, $ we have
$ (ac, ad) \in \sigma _{2}. $ Also $ (d, c) \in \sigma _{2}, $ so
$ (bd, bc) \in \sigma _{2}, $ and then $ (ac + bd, ad + bc) \in \sigma _{2}, $
i.e., $ (a, b) \ast (c, d) \in \sigma _{2}, $ also
$ (a, b) \ast (c, d) \in \sigma _{1}, $ so $ (a, b) \ast (c, d) \in
\sigma _{1} \cap \sigma _{2} \subset \rho , $ which implies that
$ (c, d) \in \rho $ because $ \rho $ is prime. So $ \sigma _{2} \subset \rho , $
i.e., $ \rho \in V^{\text{co}}(\sigma _{2}). $ This shows that
$ V^{\text{co}}(\sigma _{1} \cap \sigma _{2}) \subset
V^{\text{co}}(\sigma _{1}) \cup  V^{\text{co}}(\sigma _{2}), $ and the
equality holds. \\
(2) \ $ \rho \in \cap _{\alpha \in \Lambda } V^{\text{co}}(\sigma _{\alpha })
\Leftrightarrow \rho \in V^{\text{co}}(\sigma _{\alpha }) \ (\forall \alpha
\in \Lambda ) \Leftrightarrow \rho \supset \sigma _{\alpha } \ (\forall \alpha
\in \Lambda ) \Leftrightarrow \rho \supset
\cup _{\alpha \in \Lambda } \sigma _{\alpha } \Leftrightarrow \rho \in
V^{\text{co}}(\sigma ) $ with $ \sigma =
(\cup _{\alpha \in \Lambda } \sigma _{\alpha })^{c}. $ \\
(3) \ Obvious. The proof is completed. \quad $ \Box $ \\
From Theorem 2.25 above, the set of $ V^{\text{co}}(\sigma) $ for all
congruences $ \sigma $ on $ A $ satisfies the axiom of closed subsets,
and give a topology on $ \text{Spec}^{c}(A), $ which is called the
Zariski topology on $ \text{Spec}^{c}(A). $ A subset of $ \text{Spec}^{c}(A) $
is equipped with the subspace topology induced from the Zariski topology on
$ \text{Spec}^{c}(A). $ In the following section, a (Zariski) topology will also
be constructed for the congruence algebraic varieties.

\par     \vskip  0.4 cm

\hspace{-0.8cm}{\bf 3. Zeros of polynomial congruence equations}

\par \vskip 0.2 cm

Let $ A $ and $ B $ be two commutative semirings with $ A \subset B. $
Let $ A[x_{1}, \cdots, x_{n}] $ be the commutative semiring of polynomials
in $ n $ variables over $ A $ (see [G, p.3]). Let
$ B^{n} = \{(b_{1}, \cdots, b_{n}): \ b_{1}, \cdots, b_{n} \in B \} $
be the affine $ n-$space over $ B. $ An element $ P \in B^{n} $ will be
called a point, and if $ P = (b_{1}, \cdots, b_{n}) $ with $ b_{i} \in B, $
then the $ b_{i} $ will be called the coordinates of $ P. $
\par  \vskip 0.2cm

{\bf Definition 3.1.} \ For the commutative semirings $ A \subset B $
as above, write $ S = A[x_{1}, \cdots, x_{n}]. $ let $ T \subset
S \times S $ be a non-empty subset, and let $ \rho $ be a
congruence on $ B. $ Then we define
$$ Z_{\rho }(T)(B) =
\{P \in B^{n} : \ (f(P), g(P)) \in \rho \ \text{for all} \
(f, g) \in T \}, $$
and call $ Z_{\rho }(T)(B) $ the $ \rho-$zero set of $ T $
in $ B^{n}. $ In particular, if $ \rho = \text{id}_{B} $ is the identity
congruence on $ B, $ then $ Z_{\rho }(T)(B) $
is the set of solutions (i.e. common zeros) of the system of polynomial
equations $ \{f = g \ (\forall \ (f, g) \in T) \} $
in $ B^{n}. $ In the following, we will call $ \{(f, g) \in \rho
\ (\forall \ (f, g) \in T) \} $ a system of polynomial
$ \rho-$equations. A subset $ Y $ of $ B^{n} $ will be called an
$ \rho-$algebraic variety over $ A $ if there exists a non-empty subset
$ T \subset S \times S $ such that $ Y = Z_{\rho }(T)(B). $
\par  \vskip 0.2cm

{\bf Remark 3.2.} \ The motivation here comes from algebraic varieties
(see [Ha, Chapter 1]). As we know, one of the main problem of algebraic
geometry is to solve polynomial equations in rings. Let $ A $ and $ B $
be two commutative rings with $ A \subset B, \
\{f_{\alpha }\}_{\alpha \in \Lambda } \subset A[x_{1}, \cdots, x_{n}], $
the solutions in $ B^{n} $ of the system of polynomial
equations $ f_{\alpha } = 0 \ (\alpha \in \Lambda ) $ is then the set
$ Z(f_{\alpha })_{\alpha \in \Lambda }(B) = \{P \in B^{n} : \ f_{\alpha }(P)
= 0 \ (\forall \alpha \in \Lambda ) \} $ (see [Ha, p.2], [K, p.10]). By using
the identity congruence $ \text{id}_{B} $ on $ B, $ this set can be rewritten as
$ Z(f_{\alpha })_{\alpha \in \Lambda }(B) = \{P \in B^{n} : \ (f_{\alpha }(P),
0) \in \text{id}_{B} \ (\forall \alpha \in \Lambda ) \}. $ So it is natural
to consider the polynomial questions in semirings, in a sense of congruence,
generalizing the identity congruence (i.e. equality), as stated in the above
Def.3.1. Note that, in a ring, the equation $ f = g $ can be always changed
to be $ f - g = 0, $ i.e., the form $ h = 0. $ The case for the system of
polynomial equations $ (f_{\alpha } = g_{\alpha })_{\alpha \in \Lambda } $
is similar. But, this is not true in a semiring, because usually there is no
subtraction in a semiring. One can also consider the polynomial congruence
equations or other related congruence equations in non-commutative semirings,
yet it may be more complicated. Here we only consider the case of
commutative semirings.
\par  \vskip 0.2cm

{\bf Lemma 3.3.} \ Let $ A, B, \rho $ and $ S = A[x_{1}, \cdots, x_{n}] $
be as in Def.3.1 above. Assume $ \rho \neq B \times B. $ Then we have \\
(1) \ $ Z_{\rho }((0, 1))(B) = \emptyset , $ and $ Z_{\rho }((f, f))(B) =
B^{n} \ (\forall f \in S). $ \\
(2) \ Let $ T_{1} $ and $ T_{2} $ be two non-empty subsets of $ S \times S, $
then \\
$ Z_{\rho }(T_{1})(B) \cup Z_{\rho }(T_{2})(B) \subset Z_{\rho }
(T_{1} \ast T_{2})(B), $ where $ T_{1} \ast T_{2} $ is defined
as Def.2.3 above. \\
(3) \ The intersection of any family of $ \rho-$algebraic varieties
is an $ \rho-$algebraic variety.
\par  \vskip 0.1cm
{\bf Proof.} \ (1) \ Obvious. \\
(2) \ Let $ P \in Z_{\rho }(T_{1})(B), $ then $ (f_{1}(P), g_{1}(P)) \in
\rho , $ also $ (g_{1}(P), f_{1}(P)) \in \rho \ (\forall (f_{1}, g_{1}) \in T_{1}). $
Now for any $ \alpha \in T_{1} \ast T_{2}, $ by definition, $ \alpha =
(f_{1}, g_{1}) \ast (f_{2}, g_{2}) $ for some $ (f_{1}, g_{1}) \in T_{1} $
and $ (f_{2}, g_{2}) \in T_{2}. $ So $ \alpha (P) = (f_{1}(P), g_{1}(P))
\ast (f_{2}(P), g_{2}(P)) = (f_{1}(P)f_{2}(P) + g_{1}(P)g_{2}(P),
f_{1}(P)g_{2}(P) + f_{2}(P)g_{1}(P)) = f_{2}(P)(f_{1}(P), g_{1}(P)) +
g_{2}(P)(g_{1}(P), f_{1}(P)) \in \rho, $ and so $ P \in
Z_{\rho }(T_{1} \ast T_{2})(B), $ which shows that $ Z_{\rho }(T_{1})(B)
\subset Z_{\rho }(T_{1} \ast T_{2})(B). $ Similarly $ Z_{\rho }(T_{2})(B)
\subset Z_{\rho }(T_{1} \ast T_{2})(B). $ \\
(3) \ If $ Y_{\alpha } = Z_{\rho }(T_{\alpha })(B) $ is any family of
$ \rho-$algebraic varieties, then $ \cap Y_{\alpha } =
Z_{\rho }(\cup T_{\alpha })(B), $ so $ \cap Y_{\alpha } $ is also an
$ \rho-$algebraic variety, and the proof is completed. \quad $ \Box $
\par  \vskip 0.2cm

{\bf Theorem 3.4.} \ Let $ A, B $ and $ S = A[x_{1}, \cdots, x_{n}] $
be as in Def.3.1 above. If $ \rho \in \text{Spec}^{c}(B) $ is a prime
congruence on $ B, $ then for the $ \rho-$algebraic varieties in
$ B^{n} $ over $ A, $ we have \\
(1) \ The union of two $ \rho-$algebraic varieties is an $ \rho-$algebraic
variety; \\
(2) \ The intersection of any family of $ \rho-$algebraic varieties
is an $ \rho-$algebraic variety; \\
(3) \ The empty set $ \emptyset $ and $ B^{n} $ are $ \rho-$algebraic
varieties.
\par  \vskip 0.1cm
{\bf Proof.} \ (1) \ If $ Y_{1} = Z_{\rho }(T_{1})(B) $ and
$ Y_{2} = Z_{\rho }(T_{2})(B) $ for some $ T_{1}, T_{2} \subset
S \times S, $ then by Lemma 3.3 above, $ Y_{1} \cup Y_{2} \subset
Z_{\rho }(T_{1} \ast T_{2})(B). $ Conversely, if $ P \in
Z_{\rho }(T_{1} \ast T_{2})(B), $ and $ P \notin Y_{1}, $
then there is an element $ (f_{1}, g_{1}) \in T_{1} $ such that
$ (f_{1}(P), g_{1}(P)) \notin \rho . $ On the other hand, for any
$ (f_{2}, g_{2}) \in T_{2}, $ we have $ (f_{1}, g_{1})
\ast (f_{2}, g_{2}) \in T_{1} \ast T_{2}, $ so $ ((f_{1}, g_{1})
\ast (f_{2}, g_{2}))(P) \in \rho , $ i.e., $ (f_{1}(P), g_{1}(P))
\ast (f_{2}(P), g_{2}(P)) \in \rho , $ which implies that
$ (f_{2}(P), g_{2}(P)) \in \rho $ since $ \rho $ is a prime
congruence, so that $ P \in Y_{2}. $ Therefore
$ Y_{1} \cup Y_{2} = Z_{\rho }(T_{1} \ast T_{2})(B) $ is an
$ \rho-$algebraic variety. \\
(2) and (3) follow from Lemma 3.3 above, and the proof is completed.
\quad $ \Box $
\par  \vskip 0.2cm

{\bf Definition 3.5.} \ Let $ A, B $ and $ S = A[x_{1}, \cdots, x_{n}] $
be as in Def.3.1 above. If $ \rho \in \text{Spec}^{c}(B) $ is a prime
congruence on $ B, $ then by Theorem 3.4 above, the set of all
$ \rho-$algebraic varieties in $ B^{n} $ over $ A $ satisfies the axiom
of closed subsets, and give a topology $ \tau _{\rho , A} $ on $ B^{n}, $
i.e., a subset $ X $ of $ B^{n} $ is open in $ \tau _{\rho , A} $
if and only if its complement $ B^{n} \setminus X $ is an
$ \rho-$algebraic variety. $ \tau _{\rho , A} $ will be called the
Zariski $ \rho-$topology on $ B^{n}. $ For such $ \rho , $ a subset of
$ B^{n} $ is equipped with the subspace topology induced from
$ \tau _{\rho , A}. $ An $ \rho-$algebraic variety $ X $ of $ B^{n} $
is irreducible if it can not be expressed as the union $ X = X_{1} \cup
X_{2} $ of two proper subsets, each one of which is closed in $ X. $
\par  \vskip 0.2cm

{\bf Definition 3.6.} \ Let $ A, B, \rho $ and $ S $ be as in Def.3.1 above.
For any subset $ Y \subset B^{n}, $ we define
$$ \rho _{B}(Y) = \{(f, g) \in S \times S : \ (f(P), g(P)) \in \rho \
\text{for all} \ P \in Y \}. $$
Obviously, $ \rho _{B}(Y) $ is a congruence on $ S. $
\par  \vskip 0.2cm

{\bf Theorem 3.7.} \ Let $ A, B, \rho $ and $ S = A[x_{1}, \cdots, x_{n}] $
be as in Def.3.1 above.  \\
(1) \ If $ T_{1} \subset T_{2} $ are subsets of $ S \times S, $ then
$ Z_{\rho }(T_{1})(B) \supset Z_{\rho }(T_{2})(B). $ \\
(2) \ If $ T \subset S \times S, $ then
$ Z_{\rho }(T)(B) = Z_{\rho }(T^{c})(B) $ and $ T \subset
\rho _{B}(Z_{\rho }(T)(B)). $ \\
(3) \ If $ Y_{1} \subset Y_{2} $ are subsets of $ B^{n}, $ then
$ \rho _{B}(Y_{1}) \supset \rho _{B}(Y_{2}). $ \\
(4) \ For any two subsets $ Y_{1}, Y_{2} $ of $ B^{n}, $ we have
$ \rho _{B}(Y_{1} \cup Y_{2}) = \rho _{B}(Y_{1}) \cap \rho _{B}(Y_{2}). $
\par  \vskip 0.1cm
{\bf Proof.} \ (1) and (3) follow easily from the definitions. \\
(2) \ $ T \subset \rho _{B}(Z_{\rho }(T)(B)) $ follows directly from
the definition. For the equality, Since $ T \subset T^{c}, $ by (1),
$ Z_{\rho }(T^{c})(B) \subset Z_{\rho }(T)(B). $
Conversely, let $ P \in Z_{\rho }(T)(B), $ then
$ (h_{1}(P), h_{2}(P)) \in \rho $ for all $ (h_{1}, h_{2}) \in T. $ For any
$ (f, g) \in T^{c}, $ by Cor.2.17 above, $ f = g $ or
for some positive integer $ n $ there exists a sequence $ f = f_{1}
\rightarrow f_{2} \rightarrow \cdots \rightarrow f_{n} = g $ in which for
each $ i \in \{1, 2, \cdots , n - 1 \} $ either $ (f_{i}, f_{i + 1})
\in T^{L} $ or $ (f_{i + 1}, f_{i}) \in T^{L}. $ If $ f = g, $ then obviously
$ (f(P), g(P)) \in \rho . $ If $ (f_{i}, f_{i + 1}) \in T^{L}, $ then
$ (f_{i}, f_{i + 1}) = (a_{i}g_{i} + h_{i}, b_{i}g_{i} + h_{i}) $ for some
$ (a_{i}, b_{i}) \in T $ and some $ g_{i}, h_{i} \in S. $ By the choice
of $ P, \ (a_{i}(P), b_{i}(P)) \in \rho , $ so
\begin{align*} &(f_{i}(P), f_{i + 1}(P)) =
(a_{i}(P)g_{i}(P) + h_{i}(P), b_{i}(P)g_{i}(P) + h_{i}(P)) \\
&= (a_{i}(P), b_{i}(P)) \cdot (g_{i}(P), g_{i}(P)) + (h_{i}(P), h_{i}(P))
\in \rho .
\end{align*}
Similarly, if $ (f_{i + 1}, f_{i}) \in T^{L}, $ then
$ (f_{i + 1}(P), f_{i}(P)) \in \rho , $ which also implies that
$ (f_{i}(P), f_{i + 1}(P)) \in \rho  $ because $ \rho $ is a congruence.
So for the above sequence, we always have
$ (f_{i}(P), f_{i + 1}(P)) \in \rho $ for all
$ i \in \{1, 2, \cdots , n - 1 \}, $ that is, \\
$ (f_{1}(P), f_{2}(P)), (f_{2}(P), f_{3}(P)), \cdots , (f_{n - 1}(P), f_{n}(P))
\in \rho , $ \\
so $ (f_{1}(P), f_{n}(P)), $ i.e., $ (f(P), g(P)) \in \rho . $ This shows that
$ (f(P), g(P)) \in \rho $ for all $ (f, g) \in T^{c}, $ hence $ P \in
Z_{\rho }(T^{c})(B), $ which implies that $ Z_{\rho }(T)(B) \subset
Z_{\rho }(T^{c})(B). $ Therefore,
$ Z_{\rho }(T)(B) = Z_{\rho }(T^{c})(B). $ \\
(4) \ By (3), $ \rho _{B}(Y_{1} \cup Y_{2}) \subset \rho _{B}(Y_{i}) \
(i = 1, 2), $ so $ \rho _{B}(Y_{1} \cup Y_{2}) \subset \rho _{B}(Y_{1})
\cap \rho _{B}(Y_{2}). $ Conversely, if $ (f, g) \in \rho _{B}(Y_{1})
\cap \rho _{B}(Y_{2}), $ then $ (f(P), g(P)) \in \rho $ for all
$ P \in Y_{i} \ (i = 1, 2), $ i.e., $ (f(P), g(P)) \in \rho $ for all
$ P \in Y_{1} \cup Y_{2}, $ so $ (f, g) \in \rho _{B}(Y_{1} \cup Y_{2}), $
which implies that $ \rho _{B}(Y_{1}) \cap \rho _{B}(Y_{2}) \subset
\rho _{B}(Y_{1} \cup Y_{2}). $ Therefore, $ \rho _{B}(Y_{1} \cup Y_{2})
= \rho _{B}(Y_{1}) \cap \rho _{B}(Y_{2}). $ The proof is completed.
\quad $ \Box $
\par  \vskip 0.2cm

{\bf Lemma 3.8.} \ Let $ A, B, \rho $ and $ S = A[x_{1}, \cdots, x_{n}] $
be as in Def.3.1 above. For $ f, g \in S, P \in B^{n} $ and positive
integer $ m, $ we have \\
$ (f, g)^{\ast m}(P) = (f(P), g(P))^{\ast m} \in
B \times B $ \ (here we write $ (f, g)(P) = (f(P), g(P)) $).
\par  \vskip 0.1cm
{\bf Proof.} \ By Prop.2.6 above,
\begin{align*} &(f, g)^{\ast m} = (\sum _{2 \mid i, \ i = 0}^{m}
\begin{pmatrix}  m \\
i \end{pmatrix} f^{m - i}
 g^{i}, \ \sum _{2 \nmid i, \ i = 1}^{m} \begin{pmatrix}  m \\
i \end{pmatrix} f^{m - i} g^{i}), \ \text{so} \\
&(f, g)^{\ast m}(P) = (\sum _{2 \mid i, \ i = 0}^{m}
\begin{pmatrix}  m \\
i \end{pmatrix} f(P)^{m - i}
 g(P)^{i}, \ \sum _{2 \nmid i, \ i = 1}^{m} \begin{pmatrix}  m \\
i \end{pmatrix} f(P)^{m - i} g(P)^{i}) \\
&= (f(P), g(P))^{\ast m}.
\end{align*}
The proof is completed. \quad $ \Box $
\par  \vskip 0.2cm

Next, we will state a version of Nullstellensatz for congruences. For this,
we need write some notations. Let $ A $ be a commutative semiring and $ \theta $
be a congruence on $ A. $ Let $ S = A[x_{1}, \cdots, x_{n}] $ be the commutative
semiring of polynomials in $ n $ variables over $ A. $ For two monomials
$ a x_{1}^{i_{1}} \cdots x_{n}^{i_{n}} $ and $ b x_{1}^{i_{1}} \cdots x_{n}^{i_{n}}
\ (a, b \in A, i_{1}, \cdots , i_{n} \in \Z_{\geq 0}), $ we denote
$$ a x_{1}^{i_{1}} \cdots x_{n}^{i_{n}} \equiv b x_{1}^{i_{1}} \cdots x_{n}^{i_{n}}
\ (\text{mod} \ \theta ) \quad \text{if} \ (a, b) \in \theta . $$
In general, for two polynomials $ f = \Sigma a_{i_{1} \cdots i_{n}}
x_{1}^{i_{1}} \cdots x_{n}^{i_{n}} $ and $ g = \Sigma b_{i_{1} \cdots i_{n}}
x_{1}^{i_{1}} \cdots x_{n}^{i_{n}} $ in $ S, $ we denote $ f \equiv g \
(\text{mod} \ \theta ) $ if all the corresponding coefficients
$ a_{i_{1} \cdots i_{n}} $ and $ b_{i_{1} \cdots i_{n}} $ satisfy
$ (a_{i_{1} \cdots i_{n}}, \ b_{i_{1} \cdots i_{n}}) \in \theta . $
If $ f \equiv g \ (\text{mod} \ \theta ), $ then obviously $ (f(P), g(P))
\in \theta $ for all $ P \in A^{n}, $ where $ A^{n} $ is the affine
$ n-$space over $ A $ as before.
\par  \vskip 0.2cm

Now for the semirings $ A, B, S = A[x_{1}, \cdots, x_{n}] $ and the congruence
$ \rho $ of $ B $ in Def.3.1 above, let $ \sigma $ be a congruence on $ S, $
we define
$$ \sqrt{\sigma /\rho } = \{(f_{1}, f_{2}) \in S \times S : f_{i} \equiv g_{i}
(\text{mod} \sqrt{\rho }) \ (i = 1, 2) \ \text{for some} \ (g_{1}, g_{2}) \in
\sqrt{\sigma } \}. $$
By definition, obviously $ \sqrt{\sigma } \subset \sqrt{\sigma /\rho }, $
and if $ \sqrt{\rho } = \text{id} $ is the identity congruence, then
$ \sqrt{\sigma /\rho } = \sqrt{\sigma }. $ In general, it is easy to see that
the relation $  \sqrt{\sigma /\rho } $ on $ S $ satisfies almost all the properties
of a congruence except the transitivity, that is, it might not be transitively
closed, so it might not be a congruence. Recall that
$ (\sqrt{\sigma /\rho })^{c} $ is the congruence on $ S $ generated by the
relation $ \sqrt{\sigma /\rho }. $
\par  \vskip 0.2cm

{\bf Proposition 3.9.} \ Let $ A, B, \rho $ and $ S = A[x_{1}, \cdots, x_{n}] $
be as in Def.3.1 above. Then for any congruence $ \sigma $ on $ S, $ we have
$ (\sqrt{\sigma /\rho })^{c} \subset (\sqrt{\rho })_{B}(Z_{\rho }(\sigma )(B)). $
\par  \vskip 0.1cm
{\bf Proof.} \ Since the right-hand side of the inclusion is a congruence,
we only need to show that $ \sqrt{\sigma /\rho }
\subset (\sqrt{\rho })_{B}(Z_{\rho }(\sigma )(B)). $
let $ (f, g) \in \sqrt{\sigma /\rho }, $ then there is an element
$ (f_{1}, g_{1}) \in \sqrt{\sigma } $ such that $ f \equiv f_{1}
(\text{mod} \sqrt{\rho }) $ and $ g \equiv g_{1}
(\text{mod} \sqrt{\rho }). $ By definition,
$ (f_{1} + h, g_{1} + h)^{\ast m} \in \sigma $ for some $ h \in S $ and some
positive integer $ m. $ For any $ P \in Z_{\rho }(\sigma )(B), \
(f_{0}(P), g_{0}(P)) \in \rho $ for all $ (f_{0}, g_{0}) \in \sigma , $ so
in particular, $ (f_{1} + h, g_{1} + h)^{\ast m}(P)
\in \rho . $ Then by Lemma 3.8 above,
$ (f_{1}(P) + h(P), g_{1}(P) + h(P))^{\ast m} =
(f_{1} + h, g_{1} + h)^{\ast m}(P) \in \rho , $ so $ (f_{1}(P), g_{1}(P))
\in \sqrt{\rho}. $ which implies that $ (f(P), g(P)) \in \sqrt{\rho} $
because $ (f(P), f_{1}(P)) \in \sqrt{\rho} $ and
$ (g(P), g_{1}(P)) \in \sqrt{\rho}. $ Hence
$ (f, g) \in (\sqrt{\rho })_{B}(Z_{\rho }(\sigma )(B)), $
and the proof is completed. \quad $ \Box $
\par  \vskip 0.2cm

The Example 3.13 below gives some indication that the equality in the
above Prop.3.9 might hold. So we have the following question
\par  \vskip 0.2cm

{\bf A question about Nullstellensatz of congruences} \ Under what conditions
can the equality $ (\sqrt{\sigma /\rho })^{c} =
(\sqrt{\rho })_{B}(Z_{\rho }(\sigma )(B)) $ hold ?
\par  \vskip 0.2cm

As pointed out by an anonymous expert, recently, motivated by tropical
geometry, several versions of the Nullstellensatz question have been
raised and answered for the tropical polynomials focusing on the
case of additively idempotent semirings (see, e.g.,[BE], [JM1]). There are
other approaches were taken on this question in the literature
(see, e.g. [Gr], [GG], [Iz], [Is], [MR]).
\par  \vskip 0.2cm

{\bf Theorem 3.10.} \ Let $ A, B, \rho $ and $ S = A[x_{1}, \cdots, x_{n}] $
be as in Def.3.1 above. If $ \rho \in \text{Spec}^{c}(B) $ is a prime
congruence on $ B. $ Then for any subset $ Y \subset B^{n}, $ we have
$ Z_{\rho }(\rho _{B}(Y))(B) = \overline{Y}, $ the closure of $ Y $ in $ B^{n} $
with the Zariski $ \rho-$topology $ \tau _{\rho , A}. $
\par  \vskip 0.1cm
{\bf Proof.} \ Let $ P \in Y, $ then $ (f(P), g(P)) \in \rho $ for all
$ (f, g) \in \rho _{B}(Y), $ so $ P \in Z_{\rho }(\rho _{B}(Y))(B). $
Thus $ Y \subset Z_{\rho }(\rho _{B}(Y))(B). $ Since
$ Z_{\rho }(\rho _{B}(Y))(B) $ is closed, we get $ \overline{Y} \subset
Z_{\rho }(\rho _{B}(Y))(B). $ On the other hand, let $ W $ be any closed
subset containing $ Y. $ Then $ W = Z_{\rho }(\sigma )(B) $ for some
congruence $ \sigma $ on $ S. $ So $ Z_{\rho }(\sigma )(B) \supset Y. $
By Theorem 3.7 above, $ \sigma \subset \rho _{B}(Z_{\rho }(\sigma )(B))
\subset \rho _{B}(Y). $ So $ W = Z_{\rho }(\sigma )(B) \supset
Z_{\rho }(\rho _{B}(Y))(B). $ Therefore, $ Z_{\rho }(\rho _{B}(Y))(B)
= \overline{Y}, $ and the proof is completed. \quad $ \Box $
\par  \vskip 0.2cm

{\bf Theorem 3.11.} \ Let $ A, B, \rho $ and $ S = A[x_{1}, \cdots, x_{n}] $
be as in Theorem 3.10 above. Let $ Y \subset B^{n} $ be an $ \rho-$algebraic
variety. If $ Y $ is irreducible, then $ \rho _{B}(Y) $ is a prime congruence
on $ S. $
\par  \vskip 0.1cm
{\bf Proof.} \ Let $ (f_{1}, g_{1}), (f_{2}, g_{2}) \in S \times S. $ If
$ (f_{1}, g_{1}) \ast (f_{2}, g_{2}) \in \rho _{B}(Y), $ then by Theorems 3.7
and 3.10 above, $ Z_{\rho }((f_{1}, g_{1}) \ast (f_{2}, g_{2}))(B)
\supset Z_{\rho }(\rho _{B}(Y))(B) = \overline{Y} = Y. $ From the proof of
Theorem 3.4(1) above, we have
$ Z_{\rho }((f_{1}, g_{1}) \ast (f_{2}, g_{2}))(B) =
Z_{\rho }((f_{1}, g_{1}))(B) \cup Z_{\rho }((f_{2}, g_{2}))(B), $
so $ Y \subset Z_{\rho }((f_{1}, g_{1}))(B) \cup Z_{\rho }((f_{2}, g_{2}))(B). $
Thus $ Y = (Y \cap Z_{\rho }((f_{1}, g_{1}))(B)) \cup
(Y \cap Z_{\rho }((f_{2}, g_{2}))(B)), $ both being closed subsets of
$ Y. $ Since $ Y $ is irreducible, we have either
$ Y = Y \cap Z_{\rho }((f_{1}, g_{1}))(B), $ in which case
$ Y \subset Z_{\rho }((f_{1}, g_{1}))(B), $ or
$ Y \subset Z_{\rho }((f_{2}, g_{2}))(B). $ So by Theorem 3.7 above,
$ (f_{1}, g_{1}) \in \rho _{B}(Z_{\rho }((f_{1}, g_{1}))(B)) \subset
\rho _{B}(Y), $ i.e., $ (f_{1}, g_{1}) \in \rho _{B}(Y), $ or
$ (f_{2}, g_{2}) \in \rho _{B}(Y), $ and so $ \rho _{B}(Y) $ is a prime
congruence. The proof is completed. \quad $ \Box $
\par  \vskip 0.2cm

Let $ A, B $ be two commutative semirings, and $ \phi $
be a homomorphism from $ A $ to $ B, $ i.e.,
$ \phi : A \rightarrow B $ is a map such that
$ \phi (0) = 0, \phi (1) = 1, $ and for all
$ a, b \in A, \phi (a + b) = \phi (a )+ \phi (b) $
and $ \phi (a b) = \phi (a ) \cdot \phi (b). $ Then the kernel
$ \text{ker}\phi = \{(a_{1}, a_{2}) \in A \times A : \ \phi (a_{1}) =
\phi (a_{2}) \} $ is a congruence on $ A, $ and $ \phi $ induces
a unique injective homomorphism, say
$ \overline{\phi } : A / \text{ker}\phi \rightarrow B $ such that
$ \phi = \overline{\phi } \circ \eta , $ where $ \eta : A \rightarrow
A / \text{ker}\phi $ is the natural surjective homomorphism. Moreover,
via such $ \phi , B $ is an $ A-$semimodule, hence an $ A-$semialgebra.
In general, for a commutative semiring $ A, $ a set $ B $ is an
$ A-$semialgebra if $ B $ is both a commutative semiring and an
$A-$semimodule such that $ a (bc) = (ab)c = b(ac) \ (\forall a \in A,
b, c \in B). $ For two $ A-$semialgebras $ B $ and $ C, $ a map
$ \phi : B \rightarrow C $ is an $ A-$semialgebra homomorphism if
$ \phi $ is both a semiring homomorphism and an $ A-$semimodule
homomorphism. We let $ \text{Hom}_{A-\text{alg}}(B, C) $ to denote the
set of all $ A-$semialgebra homomorphisms from $ B $ to $ C. $ \\
Now come back to our semirings $ A \subset B $ and
$ S = A[x_{1}, \cdots, x_{n}] $ as in Def.3.1 above. Let $ \rho $ be
a congruence on $ B, $ and $ T \subset S \times S $ be a non-empty
subset. Recall that $ T^{c} $ is the congruence on $ S $
generated by $ T. $ Then it is easy to see that
both the quotient semirings $ S / T^{c} $ and $ B / \rho $ are
$ A-$semialgebras. Recall that for a subset $ Y $ of $ B^{n}, \
Y / \rho = \{(\overline{c_{1}}, \cdots , \overline{c_{n}}) : \
(c_{1}, \cdots , c_{n})  \in Y \} \subset B^{n} / \rho  $ with
$ \overline{c_{i}} = c_{i} \ \text{mod} \rho \in B / \rho . $
\par  \vskip 0.2cm

{\bf Theorem 3.12.} \ Let $ A, B $ and $ S = A[x_{1}, \cdots, x_{n}] $
be as in Def.3.1 above. Let $ \rho $ be a congruence on $ B, $ and
$ T \subset S \times S $ be a non-empty subset. Then there exists an
one-to-one map of $ Z_{\rho }(T)(B) / \rho  $ onto
$ \text{Hom}_{A-\text{alg}}(S / T^{c}, B / \rho ), $ so the cardinals
$ \sharp Z_{\rho }(T)(B) / \rho =
\sharp \text{Hom}_{A-\text{alg}}(S / T^{c}, B / \rho ). $ In particular,
if $ \rho = \text{id}_{B}, $ then $ \sharp Z_{\text{id}_{B}}(T)(B) =
\sharp \text{Hom}_{A-\text{alg}}(S / T^{c}, B). $
\par  \vskip 0.1cm
{\bf Proof.} \ By the composition of homomorphisms
$ A \hookrightarrow B \rightarrow B / \rho $ (resp.
$ A \hookrightarrow S \rightarrow S / T^{c}), \ B / \rho $ (resp.
$ S / T^{c} $) becomes a natural $ A-$semialgebra. Define a map \\
$ \phi _{0}: \ Z_{\rho }(T)(B) \longrightarrow
\text{Hom}_{A-\text{alg}}(S / T^{c}, B / \rho), \ P \mapsto \phi _{0}(P), $ \\
where $ \phi _{0}(P): \ S / T^{c} \longrightarrow B / \rho $ is defined
as follows: \\
Write $ P = (b_{1}, \cdots , b_{n}) $ with $ b_{i} \in B \ (i =1, \cdots , n), $
we have a homomorphism of semirings \\
$ \gamma _{P}: \ S \longrightarrow B / \rho , \ h(x_{1}, \cdots , x_{n})
\mapsto \overline{h(b_{1}, \cdots , b_{n})} = \overline{h(P)} \ (\forall
h \in S). $ \\
For any $ (f, g) \in T, (f(P), g(P)) \in \rho $ since $ P \in  Z_{\rho }(T)(B). $
So $ \gamma _{P}(f) = \overline{f(P)} = \overline{g(P)} = \gamma _{P}(g), $
i.e., $ (f, g) \in \text{ker} \gamma _{P}, $ so $ T \subset
\text{ker} \gamma _{P}, $ and so $ T^{c} \subset \text{ker} \gamma _{P}. $
Therefore, there is a unique homomorphism of semirings, say
$ \phi _{0}(P): \ S / T^{c} \longrightarrow B / \rho $ such that
$ \phi _{0}(P) \circ \eta = \gamma _{P}, $ where $ \eta : S \longrightarrow
S / T^{c} $ is the natural surjective homomorphism. Moreover, for any
$ a \in A $ and $ f \in S, \ \phi _{0}(P)(a \cdot \overline{f}) =
\phi _{0}(P)(\overline{af}) = \gamma _{P}(af) = \overline{(af)(P)} =
\overline{a \cdot f(P)} = \overline{a} \cdot \overline{f(P)} =
a \cdot \gamma _{P}(f) = a \cdot \phi _{0}(P)(\overline{f}), $ so
$ \phi _{0}(P) $ is also an $ A-$semimodule homomorphism, hence
$ \phi _{0}(P) \in \text{Hom}_{A-\text{alg}}(S / T^{c}, B / \rho). $
By this way, the map $ \phi _{0} $ is given. Now let
$ P, Q \in Z_{\rho }(T)(B) $ be two points such that
$ \phi _{0}(P) = \phi _{0}(Q), $ then for any
$ f \in S, \gamma _{P}(f) = \phi _{0}(P)(\eta (f)) =
\phi _{0}(P)(\overline{f}) = \phi _{0}(Q)(\overline{f}) = \gamma _{Q}(f). $
So $ \gamma _{P} = \gamma _{Q}. $ If we write $ P = (b_{1}, \cdots , b_{n}) $
and $ Q = (c_{1}, \cdots , c_{n}), $ then coordinates $ b_{i} = x_{i}(P) $
and $ c_{i} = x_{i}(Q). $ So $ \overline{b_{i}} = \overline{x_{i}(P)}
= \gamma _{P}(x_{i}) = \gamma _{Q}(x_{i}) = \overline{x_{i}(Q)} =
\overline{c_{i}}, $ and so $ \overline{P} =
(\overline{b_{1}}, \cdots , \overline{b_{n}}) = (\overline{c_{1}},
\cdots , \overline{c_{n}}) = \overline{Q}, $ i.e.,
$  \overline{P} = \overline{Q} \in Z_{\rho }(T)(B) / \rho . $
Conversely, if $ P, Q \in Z_{\rho }(T)(B) $ satisfy $ \overline{P} =
\overline{Q} \in Z_{\rho }(T)(B) / \rho , $ then obviously,
$ \gamma _{P} = \gamma _{Q}, $ and so $ \phi _{0}(P) = \phi _{0}(Q). $
Therefore, $ \phi _{0}(P) = \phi _{0}(Q) \Leftrightarrow
\overline{P} = \overline{Q}, $ and so $ \phi _{0} $ induces an injective map \\
$ \phi : \ Z_{\rho }(T)(B) / \rho \longrightarrow
\text{Hom}_{A-\text{alg}}(S / T^{c}, B / \rho), \ \overline{P} \mapsto
\phi _{0}(P) \ (\forall \ P \in Z_{\rho }(T)(B)). $ \\
Next, we define a map \\
$ \psi : \ \text{Hom}_{A-\text{alg}}(S / T^{c}, B / \rho ) \longrightarrow
Z_{\rho }(T)(B) / \rho , \ \gamma \mapsto \psi (\gamma ), $ \\
where the point $ \psi (\gamma ) \in Z_{\rho }(T)(B) / \rho $ is
defined as follows: \\
By composing $ \gamma $ with the natural $ A-$semialgebra homomorphism
$ \eta : S \longrightarrow S / T^{c}, $ we get an $ A-$semialgebra
homomorphism $ \beta = \gamma \circ \eta : \ S \longrightarrow B / \rho. $
For each $ i = 1, \cdots , n $ write $ u_{i} = \beta (x_{i}) =
\gamma (\eta (x_{i})) \in B / \rho , $ so $ u_{i} = \overline{b_{i}} $
for some $ b_{i} \in B. $ Then we define $ \psi (\gamma ) =
(\overline{b_{1}}, \cdots , \overline{b_{n}}) \in B^{n} / \rho . $
We need to show that $ (b_{1}, \cdots , b_{n}) \in
Z_{\rho }(T)(B) / \rho . $ For this, write $ P = (b_{1}, \cdots , b_{n})
\in B^{n}. $ For any $ (f, g) \in T (\subset T^{c} ), \overline{f} =
\overline{g} \in S / T^{c}, $ so $ \beta (f) = \gamma (\eta (f)) =
\gamma (\overline{f}) = \gamma (\overline{g}) = \gamma (\eta (g))
= \beta (g). $ Note that $ f = f(x_{1}, \cdots , x_{n}), g =
g(x_{1}, \cdots , x_{n}) $ and $ \beta $ is an $ A-$semialgebra
homomorphism, we have
\begin{align*} &f(\beta (x_{1}), \cdots , \beta (x_{n})) =
\beta (f(x_{1}, \cdots , x_{n})) = \beta (f) \\
&= \beta (g) = \beta (g(x_{1}, \cdots , x_{n})) =
g(\beta (x_{1}), \cdots , \beta (x_{n})),
\end{align*}
i.e., $ f(\overline{b_{1}}, \cdots , \overline{b_{n}}) =
g(\overline{b_{1}}, \cdots , \overline{b_{n}}) \in B /\rho , $
which means $ \overline{f(P)} = \overline{g(P)}, $ so
$ (f(P), g(P)) \in \rho , $ which implies that $ P \in Z_{\rho }(T)(B). $
So the above $ \psi (\gamma ) = \overline{P} \in Z_{\rho }(T)(B) / \rho . $
By this way, the map $ \psi $ is given. \\
Now we consider the composition map \\
$ \psi \circ \phi : \ Z_{\rho }(T)(B) / \rho \longrightarrow
\text{Hom}_{A-\text{alg}}(S / T^{c}, B / \rho ) \longrightarrow
Z_{\rho }(T)(B) / \rho . $ \\
For any $ P \in Z_{\rho }(T)(B), $ by definition,
\begin{align*} &(\psi \circ \phi )(\overline{P}) =
\psi (\phi (\overline{P})) = \psi (\phi _{0}(P)) =
(\phi _{0}(P)(\overline{x_{1}}), \cdots , \phi _{0}(P)(\overline{x_{n}})) \\
&= (\gamma _{P}(x_{1}), \cdots , \gamma _{P}(x_{n})) =
(\overline{x_{1}(P)}, \cdots , \overline{x_{n}(P)}) = \overline{P}.
\end{align*}
So $ \psi \circ \phi = \text{id} $ is the identity map. \\
Also for the composition map \\
$ \phi \circ \psi : \ \text{Hom}_{A-\text{alg}}(S / T^{c}, B / \rho )
\longrightarrow Z_{\rho }(T)(B) / \rho \longrightarrow
 \text{Hom}_{A-\text{alg}}(S / T^{c}, B / \rho ). $ \\
For any $ \gamma \in \text{Hom}_{A-\text{alg}}(S / T^{c}, B / \rho ), $
write $ \psi (\gamma ) = \overline{P} $ for some $ P = (b_{1}, \cdots , b_{n})
\in Z_{\rho }(T)(B). $ Then by definition, $ \psi (\gamma ) =
(\gamma (\overline{x_{1}}), \cdots , \gamma(\overline{x_{n}})), $ so
$ \gamma (\overline{x_{i}}) = \overline{b_{i}} \ (i = 1, \cdots , n). $
On the other hand, by definition, $ \phi _{0}(P)(\overline{x_{i}})
= \gamma _{P}(x_{i}) = \overline{x_{i}(P)} = \overline{b_{i}} \
(i = 1, \cdots , n). $ So $ \phi _{0}(P)(\overline{x_{i}}) =
\gamma (\overline{x_{i}}), $ and hence $ \phi _{0}(P) = \gamma , $
i.e., $ (\phi \circ \psi )(\gamma ) = \phi (\overline{P})
= \phi _{0}(P) = \gamma , $ which implies that $ \phi \circ \psi = \text{id} $
is the identity map. Therefore, both $ \phi $ and $ \psi $ are bijective,
so $ \sharp Z_{\rho }(T)(B) / \rho =
\sharp \text{Hom}_{A-\text{alg}}(S / T^{c}, B / \rho), $ and the proof is
completed. \quad $ \Box $
\par  \vskip 0.2cm

{\bf Example 3.13.} \ (1) \ Let $ A =(\Z_{\geq 0}, +, \cdot) $ be the semidomain
as before, $ p $ be a prime number, and $ \rho = \text{mod} \ p $
be the modulo $ p $ congruence, i.e.,
$ \rho = \{(a, b) : a, b \in A, a \equiv b \ (\text{mod} \ p) \}. $ Let $ S =A[t] $
be the commutative semiring of polynomials in one variable over $ A. $ Take
$ T = \{(t, 0) \} \subset S \times S, $ and let $ \sigma = T^{c} $ be the
congruence on $ S $ generated by $ T. $ By definition, we have easily that
\begin{align*} &Z_{\rho }(\sigma )(A) = Z_{\rho }(T)(A) =
\{mp : \ m \in \Z_{\geq 0}\}, \quad \text{and} \\
&(\sqrt{\rho })_{A}(Z_{\rho }(\sigma )(A)) = \{(f, g) \in S \times S :
f(0) \equiv g(0) \ (\text{mod} \ p) \}.
\end{align*}
So for $ (f, g) \in (\sqrt{\rho })_{A}(Z_{\rho }(\sigma )(A)), $ we
have $ f(0) \equiv g(0) \ (\text{mod} \ p). $ We may as well assume that
$ g(0) = f(0) + mp $ for some $ m \in \Z_{\geq 0}. $ Then
$ (f, g) = (t f_{1} + f(0), t g_{1} + f(0) + mp) $ for some $ f_{1},
g_{1} \in S. $ By Prop.2.21 above, we have
$ (t f_{1} + f(0), t g_{1} + f(0)) \in R(t,0) \subset \sigma \subset \sqrt{\sigma }. $
Note that $ f = t f_{1} + f(0) $ and
$ g \equiv t g_{1} + f(0)(\text{mod} \ p), $ in other words,
$ f \equiv t f_{1} + f(0)(\text{mod} \rho ) $ and
$ g \equiv t g_{1} + f(0)(\text{mod} \ \rho ), $ so by definition,
$ (f, g) \in (\sqrt{\sigma /\rho })^{c}, $ and so
$ (\sqrt{\rho })_{A}(Z_{\rho }(\sigma )(A)) \subset (\sqrt{\sigma /\rho })^{c}. $
Therefore, by Prop.3.9 above, the equality holds, i.e.,
$ (\sqrt{\sigma /\rho })^{c} = (\sqrt{\rho })_{A}(Z_{\rho }(\sigma )(A)). $ \\
(2) \ As discussed in (1) above, $ Z_{\rho }(T)(A) =
\{mp : \ m \in \Z_{\geq 0}\}, $ so $ Z_{\rho }(T)(A) / \rho = \{ \overline{0} \}. $
On the other hand, $ S / T^{c} = S / \sigma = A = \Z_{\geq 0} $ because
$ (t, 0) \in \sigma . $ Note that $ A / \rho =
\{\overline{0}, \cdots , \overline{p-1}\}, $ so
$ \text{Hom}_{A-\text{alg}}(S / T^{c}, A / \rho )
 = \text{Hom}_{A-\text{alg}}(A, A / \rho ) = \{\phi \}, $ where
 $ \phi : A \rightarrow A / \rho , 1 \mapsto \overline{1}. $ Hence
$ \sharp Z_{\rho }(T)(A) / \rho  = 1 =
\sharp \text{Hom}_{A-\text{alg}}(S / T^{c}, A / \rho ), $
the same as shown in Theorem 3.12 above. \quad $ \Box $
\par  \vskip 0.2cm

{\bf Remark 3.14.} \ This paper is a revised version of the early one [Q].
I thank the anonymous experts for pointing out another useful equivalent form of the
radical of congruence in Definition 2.9 above, which help me to simplify the early
version [Q, Prop.2.12, 2.13] of the results and proofs about
radical and nilpotent congruences. I also thank the anonymous experts for their useful
list of many related papers on tropical geometry
and universal algebraic geometry, that might stimulate me in further study
concerning with this work.

\par  \vskip 0.3cm

{ \bf Acknowledgments } \ I would like to thank the referee for helpful suggestions and
comments.

\par  \vskip 0.4 cm

\hspace{-0.8cm} {\bf References }
\begin{description}

\item[[BE]] A. Bertram, R. Easton, The Tropical Nullstellensatz for
congruences, Advances in Mathematics, 308 (2017), 36-82.

\item[[DMR1]] E.Yu. Daniyarova, A.G. Myasnikov, V.N. Remeslennikov,
Algebraic geometry over algebraic systems. II. Foundations. J. Math.Sci.,
185(2012), no.3, 389-416.

\item[[DMR2]] E.Yu. Daniyarova, A.G. Myasnikov, V.N. Remeslennikov,
Algebraic geometry over algebraic systems. V. The case of arbitrary signature,
Algebra Logic, 51(2012), no.1, 28-40.

\item[[G]] J.S.Golan, Semirings and Affine Equations over Them:
Theory and Applications. Boston: Kluwer Academic Publishers, 2003.

\item[[GG]] J. Giansiracusa, N. Giansiracusa, Equations of Tropical varieties,
Duke Math.J., Vol.165, 18(2016), 3379-3433.

\item[[Gr]] D. Grigoriev, On a tropical dull Nullstellensatz, arXiv: 1108.0519.

\item[[Ha]] R.Hartshorne, Algebraic Geometry. New York: Springer-Verlag,
1977.

\item[[Ho]] J.M.Howie, An Introduction to Semigroup Theory.
New York: Academic Press, 1976.

\item[[Iz]] Z. Izhakian, Tropical algebraic sets, ideals and an algebraic
Nullstellensatz, Internat.J. Algebra Comput., 18 (2008), 1067-1098.

\item[[IzS]] Z. Izhakian, E. Shustin, A tropical Nullstellensatz,
Proc. Amer. Math. Soc., 135(2007), no.12, 3815-3821.

\item[[JM1]] D. Jo\'{o}, K. Mincheva, Prime congruences of idempotent semirings
and a Nullstellensatz for tropical polynomials, arXiv: 1408.3817.

\item[[JM2]] D. Jo\'{o}, K. Mincheva, On the dimension of polynomial
semirings, arXiv: 1501.02493.

\item[[K]] G.Kemper, A Course in Commutative Algebra. Berlin:
Springer-Verlag, 2011.

\item[[MR]] D. Maclagan, F. Rinc\'{o}n, Tropical Ideals, arXiv: 1609.03838.

\item[[P1]] B. Plotkin, Some notions of algebraic geometry in universal algebra,
Algebra and Analysis 9(4)(1997), 224-248, St.Petersburg Math.j. 9(4)(1998),
859-879.

\item[[P2]] B. Plotkin,  algebraic geometry in first order logic,
Sovremennaja Matematika and Applications 22(2004), 16-62, J. Math. Sci. 137(5)
2006, 5049-5097.

\item[[PP]] B. Plotkin, E. Plotkin, Multi-sorted logic and logical geometry:
some problems. Demonstratio Mathematica, Vol.XLVIII no.4, 2015, 578-616.

\item[[Q]] Derong Qiu, On algebraic congruences, first version, arXiv: 1512.08088 v1,
26 Dec., 2015; second version, arXiv: 1512.08088 v2, 24 Jan., 2017.

\end{description}

\end{document}